\documentclass{svjour3}                     % onecolumn (standard format)
\usepackage{graphicx}
\usepackage{pgf,tikz,pgfplots}
\usepackage{tkz-euclide}
\usetikzlibrary{quotes,angles}
\usepackage{amssymb,amsmath}
\begin{document}
	\newtheorem{algo}{Algorithm}%[section]

% \[\mathcal{ABCDEFGHIJKLMNOPQRSTUVWXYZ}\]
%  \[\CMcal{ABCDEFGHIJKLMNOPQRSTUVWXYZ}\]

\title{Solving a new type of quadratic optimization problem having a joint numerical range constraint}


\titlerunning{A new type of quadratic optimization problem}        % if too long for running head

\author{Huu-Quang Nguyen       \and
        Ruey-Lin Sheu \and Yong Xia %etc.
}


\institute{Huu-Quang Nguyen \at
              Institute of Natural Science Education, Vinh
              University, Vinh, Nghe An, Vietnam\\Department of Mathematics, National Cheng Kung University, Tainan, Taiwan\\
                           \email{quangdhv@gmail.com}
           \and
           Ruey-Lin Sheu \at
              Department of Mathematics, National Cheng Kung University, Tainan, Taiwan\\
              \email{rsheu@mail.ncku.edu.tw}
              \and
              Yong Xia \at
              LMIB of the Ministry of Education, School of Mathematics and System Sciences, Beihang University, Beijing 100191, China\\
              \email{dearyxia@gmail.com}
}

\date{Received: date / Accepted: date}

\maketitle

\begin{abstract}
We propose a new formulation about quadratic optimization problems.
The objective function $F(f(x),g(x))$ is given as composition of a quadratic
function $F(z)$ with two $n$-variate quadratic functions $z_1=f(x)$ and $z_2=g(x).$ In addition, it incorporates with a set of linear inequality constraints in $z=(z_1,z_2)^T,$ while having an implicit constraint that $z$ belongs to the joint numerical range of $(f,g).$ The formulation is very general in the sense that it covers quadratic programming with a single quadratic constraint of all types, including the inequality-type, the equality-type, and the interval-type.
Even more, the composition of ``quadratic with quadratics'' as well as the joint numerical range constraint all together allow us to formulate existing unsolved (or not solved efficiently) problems into the new model. In this paper, we solve the quadratic hypersurfaces intersection problem (QSIC) proposed by P$\acute{{\rm o}}$lik and Terlaky; and the problem (AQP) to minimize the absolute value of a quadratic function over a quadratic constraint proposed by Ye and Zhang.
We show that, when $F(z)$ and the joint numerical range constraint are both convex, the optimal value of the convex optimization problem can be obtained by solving an SDP followed from a new development of the $\mathcal{S}$-procedure. The optimal solution can be approximated by conducting a bisection method on $[0,2\pi].$ On the other hand, if the joint numerical range of $f(x)$ and $g(x)$ is non-convex, the respective quadratic matrices of $f(x)$ and $g(x)$ must be linearly dependent. The linear dependence property enables us to solve (QSIC) and (AQP) accordingly by elementary analysis.

\end{abstract}
%\begin{keywords}
%  $\mathcal{S}$-procedure, Separation property, S-lemma with equality, Slater
%  condition, Intermediate value theorem, Control theory.
%\end{keywords}
%
%% REQUIRED
%\begin{AMS}
%  90C20, 90C22, 90C26, 34D20, 26B25
%\end{AMS}

\section{Introduction}
The minimization problem considered in this paper takes the following form:
\begin{eqnarray*}
{\rm (Po4)}\hspace*{2cm}
\begin{array}{lll}
v{\rm (Po4)}= &\inf\limits_{(x, z)\in \mathbb{R}^n\times \mathbb{R}^2} & F(z)  \\
 & ~~~~~{\rm s.t.} & z_1 a + z_2 b-c\leq 0,\\
& & z\in {\bf C}=\{(f(x), g(x))|~ x\in \mathbb{R}^n\},
\end{array}
\end{eqnarray*}
where $a, b, c\in \mathbb{R}^{m};$ $f(x)=x^TPx+p^Tx+p_0, g(x)=x^TQx+q^Tx+q_0$ are real-valued quadratic functions of $n$-variables with symmetric matrices $P, Q;$  $p,q\in\mathbb{R}^n$ and $p_0,q_0\in \mathbb{R}.$
The objective function is a two-variate polynomial in $z=(z_1,z_2)^T:$
\begin{eqnarray}
F(z)&=&z^T\Theta z+\eta^Tz \label{F-function}\\
&=&\theta_1z_1^2+2\theta_2z_1z_2+\theta_3z_2^2+\eta_1z_1 +\eta_2z_2, \nonumber
\end{eqnarray}
with the coefficients $\theta_1, \theta_2, \theta_3, \eta_1, \eta_2\in \mathbb{R}.$

The model {\rm (Po4)} has a very special implicit constraint set
\begin{equation}\label{C-set}
{\bf C}=\{(f(x), g(x))|~ x\in \mathbb{R}^n\}.
\end{equation} It is the joint numerical range of two quadratics $f$ and $g,$ which, in general, is a non-convex set. Even if ${\bf C}$ is convex, since it is not described specifically by convex functions (note that both $f,~g$ could be non-convex), the problem {\rm (Po4)} is still not easy to tackle. One of the major contributions in this paper is to show that, when $F$ and the joint numerical range constraint ${\bf C}$ are both convex, the (abstract)\footnote{Some researchers use the term ``abstract'' to describe the convex problem which minimizes a convex function over a convex set, while the convex set might not be described by convex functions with inequalities. The case we face here falls into the category.} convex optimization problem can be solved by an SDP followed from a new development of the $\mathcal{S}$-procedure. In this case, {\rm (Po4)} belongs to the hidden convex optimization \cite{X20}.

The model {\rm (Po4)} also has a nice extensibility so that it can be tailored to describe many existing quadratic optimization problems and creates more new ones. By writing the implicit constraint $z\in {\bf C}$ explicitly into a kind of quadratically constrained quadratic problems (abbreviated as (QCQP) in literature) in $\mathbb{R}^{n+2}:$
\begin{eqnarray}\label{nls:0}
\begin{array}{lll}
v{\rm (Po4)}= &\inf\limits_{(x, z)\in \mathbb{R}^n\times \mathbb{R}^2} & F(z)  \\
 & ~~~~~{\rm s.t.} & f(x) a + g(x) b -c\leq 0; \\
& & f(x)-z_1=0;  \\
& & g(x)-z_2=0,
\end{array}
\end{eqnarray}
we see, in the following examples, that $a f(x)+b g(x)-c\leq 0$ can be used to model quadratic inequality/equality constraints with suitable choices of $a,~b,~c,$ while $f(x)-z_1=0,~g(x)-z_2=0$ used for perturbations of the values in $f(x)$ and $g(x)$ to give meaningful applications. Interesting examples include
\begin{itemize}
	\item[(i)] When $\theta_1=\theta_2=\theta_3=0, \eta_1=1, \eta_2=0, a=0, b=1, c=0,$ the model (Po4) is reduced to a quadratic program with a single quadratic inequality constraint (abbreviated as (QP1QC) in literature):
\begin{eqnarray*}
{\rm (QP1QC)}\hspace*{2cm}
\begin{array}{lll}
&\inf\limits_{x\in \mathbb{R}^n} & f(x)  \\
 & ~~{\rm s.t.} & g(x)\leq 0,
\end{array}
\end{eqnarray*}
% p3 reviewer2
%	\end{eqnarray*}
which has been well studied in \cite{MJJ,PT07,Y03}.
	\item[(ii)] When $\theta_1=\theta_2=\theta_3=0, \eta_1=1, \eta_2=0, a=(0,0)^T, b=(1, -1)^T, c=(0, 0)^T,$ (Po4) becomes a quadratic program with a single quadratic equality constraint (abbreviated as (QP1EQC) in literature): % p2 reviewer2
	\begin{eqnarray*}%\label{po4}
		{\rm (QP1EQC)}\hspace*{0.3cm}
		\begin{array}{lll}
			&\inf\limits_{x\in \mathbb{R}^n} &f(x)\\
			&~~{\rm s.t.}&g(x)= 0.
		\end{array}
	\end{eqnarray*}
	See \cite{PW,XWS} for a complete solution to (QP1EQC).
	\item[(iii)] When $\theta_1=\theta_2=\theta_3=0, \eta_1=1, \eta_2=0, a=(0,0)^T, b=(1, -1)^T, c=(\beta, -\alpha)^T$ and $m=2$, (Po4) becomes generalized trust region subproblem (abbreviated as (GTRS) in literature): % p3 reviewer2
	\begin{eqnarray*}%\label{po4}
		{\rm (GTRS)}\hspace*{0.3cm}
		\begin{array}{lll}
			&\inf\limits_{x\in \mathbb{R}^n} &f(x)\\
			&~~{\rm s.t.}&\alpha \leq g(x) \leq\beta.
		\end{array}
	\end{eqnarray*}
	For discussions on (GTRS), please refer to \cite{MJJ,PW,WX}.
	\item[(iv)] The double well potential problems (DWP) in \cite{FGLSX,XSFX}:
	\begin{equation*}
	{\rm(DWP)} ~~\underset{x\in\mathbb{R}^n}{\min}\frac{1}{2}\left(\frac{1}{2}\|Bx-c\|^2-d \right)^2 +\dfrac{1}{2}x^TAx-u^Tx
	\end{equation*}
	can be represented by (Po4), too. Just choose
	$\theta_1=1, \theta_2=\theta_3=0, \eta_1=0, \eta_2=1, a= b=c=0,$ and let
	$f(x)=\dfrac{1}{2\sqrt{2}}\|Bx-c\|^2,$ $g(x)=-\frac{d}{2}\|Bx-c\|^2+\dfrac{1}{2}x^TAx-u^Tx+\dfrac{d^2}{2}.$
	\item[(v)] (Not solved before) In computer graphics, it is fundamental to check whether two quadratic surfaces in 3D Euclidean space intersect and to compute the intersection curve if they do. See \cite{L79,Wang-Joe-Goldman,Wilf-Manor}. P$\acute{{\rm o}}$lik and Terlaky \cite{PT07} then posed the general quadratic intersections problem in $\mathbb{R}^n.$ They ask
	``{\it Let $f(x)=x^TPx+p^Tx+p_0=0, g(x)=x^TQx+q^Tx+q_0=0$ be two quadric hypersurfaces. Can we determine whether or not the two hypersurfaces intersect without actually computing the intersections?}'' In fact, we may further assume that the input data $(P,p,p_0)$ and $(Q,q,q_0)$ are contaminated by various noises and model the general quadratic intersections problem by the following nonlinear least squares model:
	\begin{eqnarray}\label{ap:0100}
	{\rm (QSIC)}\hspace*{0.3cm}
	\begin{array}{lll}
	&\inf\limits_{(x, z)\in \mathbb{R}^{n}\times \mathbb{R}^{2}} & (z_1)^2+(z_2)^2\\
	&\hspace*{1cm}{\rm s.t.}& \begin{cases}\begin{array}{ll}
	x^TPx+p^Tx+p_0-z_1=0, \label{nls:001} \\
	x^TQx+q^Tx+q_0-z_2=0. \label{nls:002}
	\end{array}\end{cases}
	\end{array}
	\end{eqnarray}
	Obviously, (QSIC) is a special form of \eqref{nls:0}
	with $a=b=c=0.$ To our best knowledge, (QSIC) has not been solved before.
	\item[(vi)] (Not solved efficiently before) Ye and Zhang in \cite{Y03} proposed to minimize the absolute value of a quadratic function over a quadratic constraint as follows:
	\begin{eqnarray*}%\label{po4}
		{\rm (AQP)}\hspace*{0.3cm}
		\begin{array}{lll}
			&\inf\limits_{x\in \mathbb{R}^n} & |x^TPx+p^Tx+p_0|\\
			&~~{\rm s.t.}& x^TQx+q^Tx+q_0\le 0.
		\end{array}
	\end{eqnarray*}
	Note that (AQP) can be equivalently modeled as %following problem:
	\begin{eqnarray*}
		&\inf\limits_{x\in \mathbb{R}^n}& (x^TPx+p^Tx+p_0)^2   \\
		&~{\rm s.t.}& x^TQx+q^Tx+q_0\leq 0,
	\end{eqnarray*} which is again a special type of (Po4) with $\theta_1=1, \theta_2=\theta_3=0$, $\eta_1=\eta_2=0, a=c=0, b=1.$

In \cite{Y03}, (AQP) was formulated as a quadratic optimization problem subject to two quadratic constraints as follows:
\begin{equation*}
\begin{array} { c l }   \min & {t} \\ ~{\rm s.t.}& { x^TQx+q^Tx+q_0\le 0,}\\ & { -t\leq x^TPx+p^Tx+p_0 \leq t},\end{array}
\end{equation*}
and it was suggested to find the optimal value $t^*$ with the bisection method.
The procedure requires to conduct the following feasibility check for (many) fixed $t\geq 0$:
\begin{eqnarray*}%\label{po4}
		{\rm (AQP\text{-}Feas)}\hspace*{0.3cm}
		\begin{array}{lll}
			&\inf\limits_{x\in \mathbb{R}^n} & x^TQx+q^Tx+q_0\\
			&~~{\rm s.t.}& -t\le x^TPx+p^Tx+p_0\le t.
		\end{array}
\end{eqnarray*}
which is a type of (GTRS) (see case (iii) above). The entire procedure in \cite{Y03},
though polynomially implementable,
is very cumbersome due to having to execute an excessive number of SDP's.
Fortunately, our new method in this paper can resolve (AQP) by at most two  SDP's without any assumption (v.s. both primal and dual Slater conditions were assumed in \cite{Y03}).
\end{itemize}

The above examples confirm that (Po4) can be widely used in modelling many quadratic optimization problems, but, to solve it in the most general case is difficult. Interestingly, we show that, under the convexity assumptions:
\begin{eqnarray} %p2 reviewer2
\text{-}\label{cond} &&\Theta =\begin{pmatrix}\theta_1 & \theta_2 \\\theta_2 & \theta_3 \end{pmatrix} \succeq 0;\\
\text{-}\label{cond0}&&\text{the joint numerical range } {\bf C}=\{(f(x), g(x))|~ x\in \mathbb{R}^n\} \text{ is convex},
\end{eqnarray}
optimal value (Po4) can be solved by solving an (SDP). Given the convexity \eqref{cond}-\eqref{cond0}, the key is to use the separation theorem for developing a new type of $\mathcal{S}$-procedure (Theorem \ref{th1} in Sect. 2).
%
% although, historically, the latter was often used when there are more than two quadratic functions involved (such as the case in (Po4) where there are $F,~f~,g$ and $af+bg-c,$ whereas the former was mostly restricted to a system of only two quadratics.} %(with $m+2$ constraints).

The $\mathcal{S}$-procedure raises the question ``{\it when a quadratic function (think it as an objective function) restricted to a set described by quadratic functions (constraint) can be non-negative?}'' Since it can transform a quadratic optimization problem equivalently to
a family of feasibility problems and relates to Lagrange duality, it has become one of the fundamental tools in control theory and optimization. To survey the results before 2006, please refer to \cite{Derinkuyu-Pinar06,PT07}. Historically, if the constraint consists of just a single quadratic function, people refer it as $\mathcal{S}$-Lemma. If there are at least two quadratics in the constraint set, people call it the $\mathcal{S}$-procedure\footnote{Some people might prefer not to distinguishing the two terms between $\mathcal{S}$-lemma and $\mathcal{S}$-procedure though.}.

In literature, only the $\mathcal{S}$-lemma was studied completely with a necessary and sufficient condition. It includes three complete versions: the classical $\mathcal{S}$-Lemma proved in 1971 by Yakubovich (\cite{Yakubovich71}); the $\mathcal{S}$-lemma with equality proved in 2016 by Xia et. al (\cite{XWS}); and the $\mathcal{S}$-Lemma with interval bounds proved in 2015 by Wang et. al (\cite{WX}).
On the other hand, a complete version of the $\mathcal{S}$-procedure for two or more quadratic functions remains open even for two quadratic forms. Our result in this paper (Theorem \ref{th1}) is an incomplete version of the $\mathcal{S}$-procedure with $m+2$ quadratic constraints under convexity assumptions $\rm(\ref{cond})$ and $\rm(\ref{cond0}),$ but it suffices to solve (QSIC) and (AQP) completely.

To apply the new $\mathcal{S}$-procedure, we need to know in advance whether the joint numerical range {\bf C} is convex. By Theorem 4.16 in \cite{FO16}, it is known that {\bf C} can be non-convex only when the two matrices $P,~Q$ are linearly dependent. We thus solve the two problems (QSIC) and (AQP) by dividing them into two cases:
%
%
%And many incomplete versions were proposed in the literature (see \cite{XWS} for more details). Polyak (1998, \cite{PBT})  succeeded in proving an incomplete version of S-procedure for two quadratic inequalities by making an additional assumption. So far,  	conditions $\rm(\ref{cond0})$ and $\rm(\ref{cond})$ (Theorem \ref{th1} in Section 2).  Then, we show that, in Section 3, \textcolor{red}{an optimal value} $v{\rm (Po4)}$ can be computed by an SDP (Theorem \ref{thm3asdf}) and \textcolor{red}{an optimal solution} of (Po4), assuming attainable, can be obtained by a kind of bisection method on $\mathbb{R}^2$ (Lemma \ref{step1}).
%
% \textcolor{red}{We also shall} emphasize on the two applications (QSIC) and (AQP). In Section 4, we solve (QSIC) and (AQP) by dividing into two cases: %p4 #2
\begin{itemize}
	\item Suppose $\{P,Q\}$ are linearly independent; {\bf C} is convex. In this case, (QSIC) and (AQP) can be viewed as a kind of (Po4) satisfying conditions $\rm(\ref{cond})$ and $\rm(\ref{cond0}).$ The complete solution procedure is stated in Section 3.
	\item If $P=t^*Q$ (or $Q=t^*P$), by elementary analysis, we show that (QSIC) can be reduced to {\rm (QP1EQC)} (see (ii) above), while (AQP) can be reduced to finding all the solutions to a KKT system.
\end{itemize}

\section{Preliminary: A new S-procedure}

Let $f(x)=x^TPx+2p^Tx+p_0,~g(x)=x^TQx+2q^Tx+q_0$ and $F: \mathbb{R}^2\to \mathbb{R}$ be defined as in \eqref{F-function}. Given $\gamma\in\mathbb{R},$ we have the following new type of $\mathcal{S}$-procedure.
%
%Let \textcolor{red}{$F: \mathbb{R}^2\to \mathbb{R}$} be defined as in \eqref{F-function} and

\begin{theorem}\label{th1}
	Under conditions $\rm(\ref{cond})$ and $\rm(\ref{cond0})$, the following two statements are equivalent:
	\begin{itemize}
		\item[]${\rm (G_1)}$ $(\forall x\in\mathbb{R}^n,z\in\mathbb{R}^2)$ $f(x)-z_1=0, g(x)-z_2=0,$  $z_1a+z_2b\le c~\Rightarrow~ F(z)-\gamma \geq 0.$ %\label{slm}
		\item[]${\rm (G_2)}$~$(\exists \alpha,\beta\in\mathbb{R},\mu\in \mathbb{R}_+^{m})$ such that
		$F(z)-\gamma +\alpha(f(x)-z_1)+\beta(g(x)-z_2)+\mu^T(z_1a+z_2b-c)\geq 0,~~ \forall (x, z)\in\mathbb{R}^n\times\mathbb{R}^{2}$. %\label{lm}
	\end{itemize}
\end{theorem}

\proof{} Note that $\rm(G_2) \Rightarrow \rm(G_1)$ is trivial. We only prove that $\rm(G_1) \Rightarrow \rm(G_2)$. By condition (\ref{cond0}), %and Theorem \ref{th001},
${\bf C}=D_1=\{z\in \mathbb{R}^2|~z_1=f(x),~z_2=g(x),~ x\in \mathbb{R}^n\}$
is convex. Obviously,~$D_2=\{z\in \mathbb{R}^2|~z_1a+z_2b\leq c\}$ is convex
so that $D_1\cap D_2$ is convex, too. By $\rm(G_1),$
$$F(z)-\gamma\geq 0~ \forall z\in D_1\cap D_2,$$
which implies that
$$(D_1\cap D_2)\cap {D_3}=\emptyset,$$ where $$D_3:=\{z\in \mathbb{R}^2|~ F(z)-\gamma<0\}.$$
By condition (\ref{cond}), $F$ is convex,  so $D_3$ is an open convex set.

Let $\{z : v^T z = \bar{\gamma}\}$, with $v=(\bar{\alpha}, \bar{\beta})^T$, separate $D_1\cap D_2$ from $D_3.$ Since $D_3$ is open, we assume, without loss the generality, that
\begin{eqnarray}
\bar{\alpha}z_1+\bar{\beta}z_2+\bar{\gamma}\geq 0,&~\forall ~ z\in D_1\cap D_2, \label{pt4}\\
\bar{\alpha}z_1+\bar{\beta}z_2+\bar{\gamma}< 0,&~ \forall ~z\in D_3. \label{pt5}
\end{eqnarray}

From (\ref{pt5}), $\bar{\alpha}z_1+\bar{\beta}z_2+\bar{\gamma}\geq 0 \Rightarrow~ F(z)-\gamma\geq 0.$  By S-lemma, there exists $t\geq 0$ such that
\begin{equation}\label{pt5a}F(z)-\gamma-t(\bar{\alpha}z_1+\bar{\beta}z_2+\bar{\gamma})\ge 0,~ \forall ~z\in \mathbb{R}^2.\end{equation}

If $t=0$, choose $\alpha=\beta=0, \mu=0$. Then, $\rm(G_2)$ holds.

If $t>0$, by (\ref{pt4}), the following system is unsolvable:
\begin{eqnarray}
&&t\bar{\alpha}z_1+t\bar{\beta}z_2+t\bar{\gamma}<0,\nonumber \\
&&z_1a+z_2b-c\leq 0,~z\in D_1. \nonumber
\end{eqnarray}
By the Farkas theorem (see \cite[Theorem 21.1]{Rxx}, \cite[Section 6.10 21.1]{Sxx}, \cite[Theorem 2.1]{PT07}), there exists $\mu\in\mathbb{R}_+^m$ such that  \begin{equation*}\label{pt5d}
t\bar{\alpha}z_1+t\bar{\beta}z_2+t\bar{\gamma}+\mu^T(z_1a+z_2b-c)\geq 0,~\forall~z\in D_1.\end{equation*}
Equivalently, there is $\mu\ge0$ such that, $\forall~x\in \mathbb{R}^n,$
\begin{equation*}%\label{pt5e}
t\bar{\alpha}f(x)+t\bar{\beta}g(x)+t\bar{\gamma}+\mu^T(f(x)a+g(x)b-c)\geq 0.\end{equation*}

Let \begin{equation} \alpha=\mu^Ta+t\bar{\alpha}, ~~ \beta=\mu^Tb+t\bar{\beta}.\label{10fff}\end{equation}
Then, $\forall~x\in \mathbb{R}^n,$ one has:
\begin{eqnarray}
&&t\bar{\alpha}f(x)+t\bar{\beta}g(x)+t\bar{\gamma}+\mu^T(f(x)a+g(x)b-c)\geq 0 \nonumber\\ %\label{vgvg}\\
&\Leftrightarrow& (\mu^Ta+t\bar{\alpha})f(x)+(\mu^Tb+t\bar{\beta})g(x) +t\bar{\gamma}-\mu^Tc\geq 0 \nonumber \\ %\label{vfvf}\\
&\Leftrightarrow& \alpha f(x) +\beta g(x)+(\mu^Ta+t\bar{\alpha}-\alpha)z_1+(\mu^Tb+t\bar{\beta}-\beta)z_2+ t\bar{\gamma}-\mu^Tc\geq 0~~(by~ \eqref{10fff}) \nonumber \\ %\label{vcvc}\\
&\Leftrightarrow& \alpha(f(x)-z_1)+\beta(g(x)-z_2)+\mu^T(z_1a+z_2b-c) \geq -t\bar{\alpha}z_1-t\bar{\beta}z_2-t\bar{\gamma}. \label{pt5f}
\end{eqnarray}
%
%Eq. \eqref{vfvf} is obtained from \eqref{vgvg} by changing the order and combining factors. Eq. \eqref{vcvc} holds true by the fact that $\mu^Ta+t\bar{\alpha}=\alpha$, $\mu^Tb+t\bar{\beta}=\beta$ and  $\mu^Ta+t\bar{\alpha}-\alpha=\mu^Tb+t\bar{\beta}-\beta=0$ (see \eqref{10fff}). Eq. \eqref{pt5f} is obtained from \eqref{vcvc} by changing the order, combining and moving some factors to the right hand of sign $\geq.$

Finally, we combine (\ref{pt5a}) with (\ref{pt5f}) to obtain $\rm(G_2)$. \hfill$\Box$
\endproof

As for condition (\ref{cond0}), there is an easy-to-verify sufficient condition as follows.

%
%Notice that, the joint numerical range {\bf C} of $\{f, g\}$ be defined as in \eqref{C-set}.
%Let us first quote a known result concerning the convexity of {\bf C}.

\begin{theorem}[Theorem 4.16 in \cite{FO16}]\label{th001}
If $\{P,Q\}$ are linearly independent, the joint numerical range ${\bf C}$ defined by \eqref{C-set} is a convex set in $\mathbb{R}^2$.
\end{theorem}

In the following, we give two examples to show that conditions (\ref{cond}) and (\ref{cond0}) cannot be omitted from Theorem \ref{th1}.

\begin{example} Let $f(x)=x_1+x_2,~g(x)=2x_1^2-x_2^2,$ $~F(z)=4z_1^2+z_2$ and $a=b=c=\gamma=0.$ Note that $F(z)$ is convex and condition (\ref{cond}) holds.
However, condition (\ref{cond0}) is violated since the joint numerical range
	\begin{equation*}
	{\bf C}=\{z\in\mathbb{R}^2~|~ z_2\geq -2z_1^2\}
	\end{equation*}
is not convex. See \cite[Example 3.1]{PBT}.

%
%Since $P=0,$ by Theorem \ref{th001},
%	
	It can be verified that  $z_1=x_1+x_2 \text{~and~} z_2=2x_1^2-x_2^2$ imply that
	$$F(z)-0=4(x_1+x_2)^2+(2x_1^2-x_2^2)=6x_1^2+8x_1x_2+3x_2^2\geq 0.$$
	Therefore, $\rm(G_1)$ holds. On the other hand,
	\begin{eqnarray*}
		&& F(z)-0+\alpha(f(x)-z_1)+\beta(g(x)-z_2) \\
		&=&4z_1^2+z_2+\alpha(x_1+x_2-z_1)+\beta(2x_1^2-x_2^2-z_2)\geq 0
	\end{eqnarray*}
	holds if and only if
	\begin{equation*}\label{inequa1}M=\begin{pmatrix}
	2\beta & 0 &~ 0 &0&\alpha/2\\
	0 & -\beta & 0 & 0&\alpha/2\\
	0 & 0 & 4 & 0&-\alpha/2\\
	0 & 0 & 0 & 0&(1-\beta)/2\\
	\alpha/2 &~~ ~~\alpha/2 & ~~-\alpha/2 & (1-\beta)/2&0
	\end{pmatrix}\succeq 0.\end{equation*}
	To make $M\succeq0,$ we need $\beta=0$.
	However, it leads to $$\begin{pmatrix}
	0&(1-\beta)/2\\
	(1-\beta)/2&0
	\end{pmatrix}=\begin{pmatrix}
	0&1/2\\
	1/2&0
	\end{pmatrix}\not\succeq0$$ and thus $\rm(G_2)$ fails. \end{example}

\begin{example} Let $f(x)=x_1^2, g(x)=x_2^2,$ and $F(z)=2z_1z_2;$
	$a=b=c=\gamma=0.$ By Theorem \ref{th001}, condition (\ref{cond0}) holds.
However, condition (\ref{cond}) is clearly violated.
	
	Note that $z_1=x_1^2$ and $z_2=x_2^2$ imply that $F(z)-0 =x_1^2x_2^2\geq 0.$ Therefore $\rm(G_1)$ holds. On the other hand,
	\begin{eqnarray*}
		&&F(z)-0+\alpha(f(x)-z_1)+\beta(g(x)-z_2) \\
		&=&2z_1z_2+\alpha(x_1^2-z_1)+\beta(x_2^2-z_2)\geq 0
	\end{eqnarray*} holds if and only if
	\begin{equation*}\label{inequa2}M=\begin{pmatrix}
	\alpha & ~0 & 0 &0&0\\
	0 &~ \beta & 0 & 0&0\\
	0 &~ ~0 & 0 & 1&\frac{-\alpha}{2}\\
	0 & ~ ~0 & 1 & 0&\frac{-\beta}{2}\\
	0 &~ ~0 & \frac{-\alpha}{2} & \frac{-\beta}{2}&0
	\end{pmatrix}\succeq 0.\end{equation*}
	However, this is impossible since
	$\begin{pmatrix}
	0&1\\
	1&0
	\end{pmatrix}\not\succeq 0.$ Then, $\rm(G_2)$ fails.
\end{example}

\section{Solving (Po4) under conditions (\ref{cond}) and (\ref{cond0})}
\subsection{Computing the optimal value $v(\rm{Po4})$}

Applying the new $\mathcal{S}$-procedure in Theorem \ref{th1}, we show that the optimal value $v(\rm{Po4})$ under conditions (\ref{cond}) and (\ref{cond0}) can be obtained by solving an SDP.

\begin{theorem}\label{thm3asdf} Under conditions $\rm(\ref{cond})$ and $\rm(\ref{cond0})$, the optimal value of $\rm (Po4)$, $v(\rm{Po4}),$ can be computed by
	\begin{align}\label{19aa}
	v(\rm{Po4})=\underset{\tiny \begin{array}{lll}
		\gamma, \,\alpha, \,\beta \in \mathbb{R}\\
		\mu\in \mathbb{R}_+^m
		\end{array}}{\sup}  \left\lbrace\begin{array}{cc} \gamma&|~M \succeq 0 \end{array}\right\rbrace,
	\end{align}
	where $M\in\mathbb{R}^{(n+3)\times(n+3)}$ is
	\begin{equation}\label{sigma.1}
	\begin{pmatrix}\begin{array}{ccc}
	\theta_1 & \quad \quad ~& \theta_2 \\ \theta_2 &\quad & \theta_3
	\end{array} & ~~[0]_{2\times n} & \begin{array}{cc}
	\frac{\mu^Ta+\eta_1-\alpha}{2}  \\ \frac{\mu^Tb+\eta_2-\beta}{2}
	\end{array}\\  ~~~[0]_{n\times 2} & \quad \alpha P+\beta Q & \alpha p+\beta q\\ \begin{array}{cc}
	\frac{\mu^Ta+\eta_1-\alpha}{2}  & \quad \frac{\mu^Tb+\eta_2-\beta}{2}
	\end{array} & \quad ~\alpha p^T+\beta q^T & \quad \alpha p_0+\beta q_0-\mu^Tc-\gamma\end{pmatrix},
	\end{equation}where
	 $[0]=(0_{ij})_{2\times n}$ with $0_{ij}=0 ~\forall i, j.$ \end{theorem}

\proof
Since (Po4) can be formulated as \eqref{nls:0}, we have
%
%
%\begin{eqnarray*}
%{\rm (Po4)}\hspace*{2cm}
%\begin{array}{lll}
%v{\rm (Po4)}= &\inf\limits_{(x^T, z_1, z_2)^T} & F(z_1,z_2)  \\
% & ~~~~~{\rm s.t.} & z_1 a + z_2 b-c\leq 0,\\
%& & (z_1,z_2)\in {\bf C}=\{(f(x), g(x)): x\in \mathbb{R}^n\}
%\end{array}
%\end{eqnarray*}

\begin{equation}
\begin{array}{lll}
v{\rm (Po4)}&=&\inf\limits_{(x, z)\in \mathbb{R}^n\times \mathbb{R}^2} F(z) \\
& ~~~{\rm s.t.} & \begin{cases}\begin{array}{lll}
f(x)a+g(x)b-c\leq 0 \\
f(x)-z_1=0 \\
g(x)-z_2=0
\end{array}\end{cases}\\
&=&{\sup} \left\lbrace \gamma \left|~ \left\lbrace
(x, z)\in \mathbb{R}^n\times \mathbb{R}^2 \left|~ \begin{array}{ll}
F(z) < \gamma\\
z_1a+z_2b\leq c\\
f(x)=z_1\\
g(x)=z_2
\end{array}\right.
\right\rbrace=\emptyset\right.\right\rbrace\\
%&=&{\sup} \left\lbrace \gamma\left|~  \left\lbrace
%\begin{array}{ll}
%z_1a+z_2b\leq c\\
%f(x)=z_1\\
%g(x)=z_2
%\end{array}
%\Rightarrow F(z_1,z_2) \geq \gamma \right\rbrace\right.\right\rbrace \\
&=&\underset{\tiny\begin{array}{lll}
	\gamma, \,\alpha, \,\beta \in \mathbb{R}\\
	\mu\in \mathbb{R}_+^m
	\end{array}}{\sup} \big\{ \gamma |~\overline{F}(x,z)\geq 0~\forall (x, z)\in \mathbb{R}^n\times \mathbb{R}^2\big\},
\end{array}\label{19uio}
\end{equation}
where $\overline{F}(x, z)=F(z)-\gamma +\alpha(f(x)-z_1)+\beta(g(x)-z_2)+\mu^T(z_1a+z_2b-c)$, and the last equality in (\ref{19uio}) holds by Theorem \ref{th1}. Note that $\overline{F}(x,z)\ge0$ can be written as a linear matrix inequality \eqref{19aa} with a matrix $M$ defined in \eqref{sigma.1}. \hfill$\Box$  \endproof

To illustrate Theorem \ref{thm3asdf}, we provide a numerical example below.

\begin{example}
	Let
	$P=\begin{pmatrix}1 & 0 & 0\\0 & 2 & 0\\0 & 0 & 3\end{pmatrix},$ $p=\begin{pmatrix}0\\1\\1\end{pmatrix},$ $p_0=7$, $Q=\begin{pmatrix}1 & -2 & 2\\-2 & 1 & 3\\2 & 3 & 1\end{pmatrix},$ $q=\begin{pmatrix}1\\2\\3\end{pmatrix},$ $q_0=2$, $\Theta=\begin{pmatrix}1 & 0\\0 & 2\end{pmatrix},$ $\begin{pmatrix}\eta_1\\ \eta_2\end{pmatrix}=\begin{pmatrix}1\\2\end{pmatrix}$,
	$a=b=c=0.$
	The SDP (\ref{19aa}) becomes:
	
	$$\underset{\gamma, \,\alpha, \,\beta \in \mathbb{R}}{\max}  \left\lbrace\begin{array}{cc} \gamma&|~M\succeq 0 \end{array}\right\rbrace,$$ where $M$ is \[\begin{pmatrix}\begin{array}{cccccc}
	~1~&~\quad 0~&\quad 0&\quad 0&\quad0&\frac{1-\alpha}{2}\\
	~0~&~\quad 2~&\quad0&\quad 0&\quad 0&\frac{2-\beta}{2}\\
	~0~&~\quad 0~&\quad\alpha+\beta&\quad-2\beta&\quad 2\beta&\beta\\
	~0~&~\quad 0~&~ -2\beta&\quad 2\alpha+\beta&\quad 3\beta&\alpha+2\beta\\
	~0~&~\quad 0~&\quad 2\beta&3\beta&3\alpha+\beta&\alpha+3\beta\\
	\frac{1-\alpha}{2}&\quad \frac{2-\beta}{2}&\quad \beta&\quad \alpha+2\beta&\quad \alpha
	+3\beta&\quad {7\alpha+2\beta-\gamma}
	\end{array} \end{pmatrix}.\]   The optimal value (cvx$_-$optval): +43.7102.\end{example}
\begin{remark} The computational experiment was conducted in Matlab version R2016a running on a PC with Core i5 CPU and 8G memory. The SDP program (\ref{19aa}) is modeled by CVX 1.21 \cite{BOY} and solved by SDPT3 within CVX.
\end{remark}

\subsection{Solving for an optimal solution of (Po4)}

In general, problem (Po4) under conditions (\ref{cond}) and (\ref{cond0}) may not have an optimal solution. First, it can be unbounded from below. For example, let $f(x)=x_1^2, g(x)=x_2^2$, $F(z)=z_1^2-z_2$, $a=b=c=0$. Then both conditions (\ref{cond}) and (\ref{cond0}) hold. It is easy to see that
\begin{equation*}v({\rm{Po4}})={\inf_{x\in \mathbb{R}^2}\left\{x_1^4-x_2^2\right\}=-\infty.}
\end{equation*}
Secondly, (Po4) may be bounded but not attainable. For example, let $f(x)=x_1^2,$ $g(x)=x_1x_2-1,$ $F(z)=z_1^2,$ $a=c=(0,0)^T, b=(1,-1)$. Then (Po4) becomes \begin{eqnarray*}
	&\inf\limits_{x\in \mathbb{R}^2}& x_1^4   \\
	&{\rm s.t.}& x_1x_2-1=0.
\end{eqnarray*}
The optimal value is zero, but there exists no $x\in \{x| x_1x_2=1\}$ such that $x_1^4=0$.

 Even if the optimal value $v(\rm{Po4})$ is attained at some $x^*,$ we emphasize that it is generally not easy to solve $x^*$ from the related SDP optimal solution. Ideally, we can use the following procedure for obtaining $x^*$.
   %It is obvious that $v{\rm (Po4)}$ is the minimal value such that $\{z\in \mathbb{R}^2: F(z)\leq v{\rm (Po4)}\}\cap \{z\in \mathbb{R}^2: z\in {\bf C}, z_1a+z_2b\leq c\}\ne \emptyset$. On the other hand,  given $\epsilon>0$, we get an approximation value of $v{\rm (Po4)}$,  ${q}_0\in [v{\rm (Po4)}-\epsilon, v{\rm (Po4)}+\epsilon]$ by solving the $\rm (SDP)$ in \eqref{19aa}. One has ${q}_0+\epsilon\geq v{\rm (Po4)}$. Therefore,  the value ${q}_0+\epsilon$ is to be guaranteed on the overall solution $x^*$ in \eqref{subp41}.

   \begin{algo}\label{algo1}  Solve an optimal solution $x^*$ to {\rm(Po4)} under conditions \eqref{cond}-\eqref{cond0}.\vskip0.2cm	
%\begin{itemize}
%  \item{\bf Step 1:} {\rm Compute the optimal value $v(\rm{Po4})$  by solving the $\rm (SDP)$ in \eqref{19aa}.}
%  \item{\bf Step 2:} {\rm Solve $(z_1^*, z_2^*)^T\in\mathbb{R}^2$ that satisfies
%	\begin{eqnarray}\label{subp4}
%	\left\{\begin{array}{ll}
%	F(z^*_1, z^*_2)=v(\rm{Po4});\\
%	z^*_1a+z^*_2b\leq c;\\
%	(z_1^*,z_2^*)\in {\bf C},
%	\end{array}\right.
%	\end{eqnarray}
%	where ${\bf C}$ is the joint numerical range defined by \eqref{C-set}.}
%  \item{\bf Step 3:} {\rm Solve $x^*\in\mathbb{R}^n$ that satisfies
%	\begin{eqnarray}\label{subp41}
%	\left\{\begin{array}{ll}
%	f(x^*)=z_1^*\\
%	g(x^*)=z_2^*.
%	\end{array}\right.
%	\end{eqnarray}}
%\end{enumerate}
\noindent {\bf Step 1:} {\rm Compute the optimal value $v(\rm{Po4})$ by solving the $\rm (SDP)$ in \eqref{19aa}.}\\
\noindent {\bf Step 2:}  {\rm Solve $z^*\in\mathbb{R}^2$ that satisfies
	\begin{eqnarray}\label{subp4}
	\left\{\begin{array}{ll}
	F(z^*)=v(\rm{Po4});\\
	z^*_1a+z^*_2b\leq c;\\
	z^*\in {\bf C},
	\end{array}\right.
	\end{eqnarray}
where ${\bf C}$ is the joint numerical range defined by \eqref{C-set}.}\\
\noindent {\bf Step 3:} {\rm Solve $x^*\in\mathbb{R}^n$ that satisfies
	\begin{eqnarray}\label{subp41}
	\left\{\begin{array}{ll}
	f(x^*)=z_1^*,\\
	g(x^*)=z_2^*.
	\end{array}\right.
	\end{eqnarray}}
\end{algo}

Note that Step 3 is to solve a simultaneous system of two quadratic equations. Theoretically, it is known that a solution to a system of $k$ quadratic equations, $k$ fixed, can be found in polynomial time.
See, e.g., the paper by Grigoriev and Pasechnik in 2005 \cite{Grigoriev-Pasechnik05}. However, their procedure appears to be difficult for implementation. In practice, we suggest to use the Newton method for finding the root $x^*$ in \eqref{subp41}.

Let us focus on Step 2 in Algorithm \ref{algo1}. Under condition (\ref{cond}), $F(z)$ is a convex function. Then,
$\{z\in \mathbb{R}^2|~ F(z)=v(\rm{Po4})\}$
is either an ellipse, a parabola, or a pair of parallel straight lines in $\mathbb{R}^2$, whose forms are either $F(z)=\theta_1 z_1^2+\theta_3z_2^2$; $F(z)=\theta_1z_1^2+\eta_2z_2, \eta_2\ne 0$ or $F(z)=\theta_1z_1^2$, respectively. Figure 1 depicts the case for
$F(z)=\theta_1 z_1^2+\theta_3z_2^2$ (in black dash), while the set (in stripes bounded by solid curves)
\begin{equation}\label{sub-3.4}
\Omega_0={\bf C}\cap \{z\in \mathbb{R}^2 \vert~ z_1a+z_2b\leq c\}
\end{equation}
is also convex under condition (\ref{cond0}). The intersection
$$z^*\in \Omega_0 \cap \{z\in \mathbb{R}^2|~ F(z)=v(\rm{Po4})\} $$
is the required solution for the subproblem (\ref{subp4}).
Below we propose a type of bisection method on $[0,2\pi]$ for computing an approximate solution $\bar z$ of $z^*.$

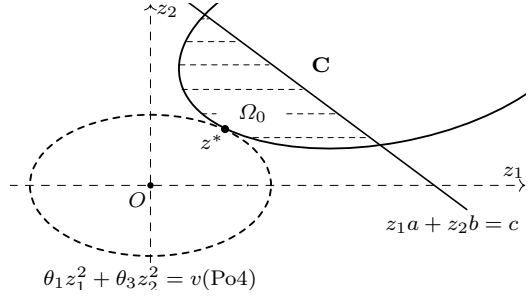
\begin{figure}   %line width=1.pt,dash
	\centering
	\begin{tikzpicture}[line cap=round,line join=round,scale=0.32]\label{fig.01}
	\clip(-5.7458,-4.33) rectangle (15.4,7.548);
	
	\draw [color=black,rotate around={0.:(0.,0.)},line width=0.7pt,dash pattern=on 2.3pt off 2.3pt] (0.,0.) ellipse (4.956054411039317cm and 2.926170761452972cm);
	\draw [color=black,rotate around={10.445622:(9.4895238,5.96658)},line width=0.7pt] (9.4895, 5.96658) ellipse (8.4242cm and 4.2379cm);
	\draw (0,-3.89) node {${\theta_1z_1^2+\theta_3 z_2^2=v(\rm{Po4})}$};
	\draw (7,5) node {$\bf{C}$};
	\draw (4.2*0.6,3*0.6) node {${z^*}$};
	\draw (4.2,3) node {${\Omega_0}$};
	\draw [->,line width=0.3pt,dash pattern=on 3pt off 3pt] (0,-3.2) -- (0,7.5);
	\draw [->,line width=0.3pt,dash pattern=on 3pt off 3pt] (-5.8,0) -- (15.4,0);
	\draw [color=black, line width=0.6pt] (13,-1) -- (-5,12.5);
	\draw (12.4,-1.6) node {${z_1a+z_2b= c}$};
	\draw (14.93,0.534) node {${z_1}$};
	\draw (0.7,7.21) node {${z_2}$};
	\draw (-0.53,-0.59) node {$O$};	

	\draw [line width=0.2pt,dash pattern=on 1pt off 1pt] (2.1, 6.9)-- (4.98-2.6, 6.9);
	\draw [line width=0.2pt,dash pattern=on 2pt off 2pt] (1.5, 5.957)-- (6.4079-2.8, 5.957);
	\draw [line width=0.2pt,dash pattern=on 2pt off 2pt] (1.2, 4.995)-- (7.673-2.8, 4.995);
	\draw [line width=0.2pt,dash pattern=on 2pt off 2pt] (1.3, 3.9766)-- (9.0312-2.8, 3.9766);
	\draw [line width=0.2pt,dash pattern=on 1.5pt off 2pt] (2.15, 2.9669)-- (5.7-2.8, 2.9669);
	\draw [line width=0.2pt,dash pattern=on 2pt off 2pt] (5.6, 2.9669)-- (10.377-2.8, 2.9669);
	\draw [line width=0.2pt,dash pattern=on 2pt off 2pt] (4.22, 1.98)-- (11.78-2.8, 1.98);
		
	\begin{scriptsize}
	\draw [color=black] (4,0) circle (0.5pt);
	\draw [fill=black] (3.064174, 2.33326) circle (4.1pt);
	\draw [fill=black] (0, 0) circle (3.1pt);
	\end{scriptsize}	
	\end{tikzpicture}
	\caption{Graphic representation for the subproblem \eqref{subp4}.}
\end{figure}

For simplicity, let us assume that $F(z)=z_1^2+z_2^2$ for Algorithm 2 below.
Given $\epsilon>0$ and suppose in Step 1 of Algorithm \ref{algo1} we obtain $\bar{v}$ as an approximate value of $v(\rm{Po4})$ by solving the SDP in \eqref{19aa} such that \begin{equation}\label{vbar}
\bar{v}\in [v({\rm{Po4}}), v({\rm{Po4}})+\epsilon/2].
\end{equation}

%After that we will find an approximation solution $\bar{x}$ of (Po4) such that $v(Po4)\leq F(f(\bar{x}), g(\bar{x}))\leq v(Po4)+\epsilon$.

%\begin{definition}\label{def1} \textcolor{red}{The point $z\in\mathbb{R}^2$ is called to be zxzxzx if $[Oz]\cap{\Omega_0}\ne\emptyset.$    }	\end{definition}

\begin{algo}\label{algo2} Implementation of Step 2 in Algorithm \ref{algo1}, assuming $F(z)=z_1^2+z_2^2.$
\vskip0.1cm	

\noindent{\bf Step 2.1:} {\rm Let $\Omega_0$ as in \eqref{sub-3.4} and $\epsilon,~\bar{v}$ as in \eqref{vbar}. Set $l_0:=0, u_0:=2\pi,~k:=0.$}
%where $\bar{v}$ satisfies \eqref{vbar}.  %$\alpha_{\bar{v}}=\min\{\arcsin\dfrac{\bar{v}}{\bar{v}+\epsilon/2}, \pi/4\}$  $p_{\bar{v}}=\left[\log_2\arcsin\dfrac{\bar{v}}{\bar{v}+\epsilon/2}\right]+1$

\noindent{\bf Step 2.2:} {\rm Set $\varphi_{k+1}:=\dfrac{l_k+u_k}{2}$ and
$$\Omega_{k+1}:=\Omega_k\cap \{z\in\mathbb{R}^2\vert ~ (\sin\varphi_{k+1}, -\cos\varphi_{k+1})z\ge 0\}.$$
Solve the following problem by formulating it as an SDP in the form of \eqref{19aa}:
\begin{equation*}
	v_{k+1}=\inf\big\{F(z)\big|~z\in \Omega_{k+1}\big\}.
\end{equation*}}
\noindent{\bf Step 2.3:} {\rm Test whether $z^*$ belongs to $\Omega_{k+1}:$}
\begin{itemize}
  \item {\rm If $v_{k+1}\le v\rm(Po4),$ set $u_{k+1}:=\varphi_{k+1},~l_{k+1}:=l_{k}.$} %,~k:=k+1$ and go to Step 2.}
  \item {\rm Otherwise, set $u_{k+1}:=u_{k},~l_{k+1}:=\varphi_{k+1}$ and $$\Omega_{k+1}:=\Omega_k\cap \{z\in\mathbb{R}^2\vert ~ (\sin\varphi_{k+1}, -\cos\varphi_{k+1})z\le 0\}.$$}  %Set $k:=k+1$ and go to Step 2.}
\end{itemize}

\noindent{\bf Step 2.4:} {\rm Set
	$k:=k+1.$ If $\vert u_{k}-l_{k}\vert > \arccos\frac{\sqrt{\bar{v}}}{\sqrt{\bar{v}+\epsilon/2}},$ go to Step 2.2. Otherwise, set
$\check{z}=(\sqrt{\bar{v}}\cos u_{k},\sqrt{\bar{v}}\sin u_{k}),$ $\hat{z}=(\sqrt{\bar{v}}\sec{\frac{2\pi}{2^{k}}}\cos l_{k},\sqrt{\bar{v}}\sec{\frac{2\pi}{2^{k}}}\sin l_{k})$.}
%$\check{z}$, where $\check{z}_1=\bar{v}\cos u_{k}, \check{z}_2= \bar{v}\sin u_{k}$, }

\noindent{\bf Step 2.5:} {\rm Test whether $[O,\check{z}]\cap{\Omega_0}\ne\emptyset$ or $[O,\hat{z}]\cap{\Omega_0}\ne\emptyset$ :}
	\begin{itemize}
\item {\rm if $[O,\check{z}]\cap{\Omega_0}\ne\emptyset$ (Fig. 2), find a point $\bar{z}$ in $\{z: \check{z}_2z_1-\check{z}_1z_2=0\}\cap \Omega_0$ nearest to $O$ (if $[O,\hat{z}]\cap{\Omega_0}\ne\emptyset$, find a point $\bar{z}$ in $\{z: \hat{z}_2z_1-\hat{z}_1z_2=0\}\cap \Omega_0$ nearest to $O$).   Report $\bar{z}$ as an approximate solution of $z^*$.}
\item {\rm  Otherwise (Fig. 3): if $[O,\check{z}]\cap{\Omega_0}=\emptyset$ and $[O,\hat{z}]\cap{\Omega_0}=\emptyset$, find a point $\bar{z}$ in $\{z: \check{z}_1z_1+\check{z}_2z_2={\bar{v}}\}\cap \Omega_0$ nearest to $\check{z}$. }
	\end{itemize}
 \end{algo}

\begin{figure}   %line width=1.pt,dash
	\centering
	\begin{tikzpicture}[line cap=round,line join=round,scale=0.42]\label{fig.1a}
	\clip(-5.7458,-5) rectangle (15.4,7.548);
	
	\draw [color=black,rotate around={0.:(0.,0.)},line width=0.6pt,dash pattern=on 2.1pt off 2.1pt] (0.,0.) ellipse (3.99cm and 3.99cm);
	
	\draw [color=black,rotate around={10.445622:(9.4895238,5.96658)},line width=0.6pt] (9.4895, 5.96658) ellipse (8.4242cm and 4.2379cm);
	
	\draw [->,line width=0.3pt,dash pattern=on 3pt off 3pt] (0,-3.2) -- (0,7.5);
	\draw [->,line width=0.3pt,dash pattern=on 3pt off 3pt] (-5.8,0) -- (15.4,0);
	\draw [color=black, line width=0.7pt] (13,-1) -- (-5,12.5);

	\draw [line width=0.2pt,dash pattern=on 1pt off 1pt] (2.1, 6.9)-- (4.98-2.6, 6.9);
	\draw [line width=0.2pt,dash pattern=on 2pt off 2pt] (1.5, 5.957)-- (6.4079-2.8, 5.957);
	\draw [line width=0.2pt,dash pattern=on 2pt off 2pt] (1.2, 4.995)-- (7.673-2.8, 4.995);
	\draw [line width=0.2pt,dash pattern=on 2pt off 2pt] (1.3, 3.9766)-- (9.0312-2.8, 3.9766);
	\draw [line width=0.2pt,dash pattern=on 2pt off 2pt] (4.9, 2.9669)-- (10.377-2.8, 2.9669);
	\draw [line width=0.2pt,dash pattern=on 2pt off 2pt] (4.22, 1.98)-- (11.78-2.8, 1.98);	
	
	\draw [line width=0.5pt,color=black,dash pattern=on 1pt off 1pt] (3.17,2.53)-- (10, 8);
	
	\draw [line width=0.5pt,color=black,dash pattern=on 1pt off 1pt] (2.13*1.043, 3.17*1.045)-- (10, 15);
		
	\begin{scriptsize}
	\draw [fill=black] (2.33, 2.7999) circle (1.6pt);
	\draw [fill=black] (0, 0) circle (1.6pt);	
	\draw [fill=black] (2.13*1.043, 3.17*1.045) circle (1.6pt); %\check z
	\draw [fill=black] (2.135*0.953, 3.17*0.956) circle (1.6pt);
	\end{scriptsize}
	
	\draw (14.93,0.35) node {${z_1}$};
	\draw (0.5,7.21) node {${z_2}$};
	
	\draw (1.96*1.1,2.2*1.1) node {${z^*}$};
	\draw (2.17,3.67) node {$\check{z}$};
	\draw (1.99*0.83,3.6*0.83) node {$\bar{z}$};
	
	\draw (3.3,3.5) node {$\Omega_0$};
	\draw (-0.5,-0.5) node {$O$};
		
	\draw (12.4,-1.6) node {${z_1a+z_2b= c}$};
	\draw (7.4,-2.8) node {where $\varphi\leq\arccos\frac{\sqrt{\bar{v}}}{\sqrt{\bar{v}+\epsilon/2}}$};
	\draw (0,-4.5) node {${z_1^2+ z_2^2=\bar{v}\in[v(Po4), v(Po4)+\epsilon/2]}$};
	
	% danh dau goc
	\draw
	(3.17,2.53) coordinate (a) node[right] {}
	-- (0,0) coordinate (b) node[left] {}
	-- (2.131*1.02,3.178*1.02) coordinate (c) node[above right] {}
	pic["$\varphi$",draw=black,<->,angle eccentricity=1.2,angle radius=0.8cm] {angle=a--b--c};
	%het danh dau goc	
	\end{tikzpicture}
	\caption{Graphic representation for the subproblem \eqref{subp4}: $[O,\check{z}]\cap{\Omega_0}\ne\emptyset$ .}
\end{figure}
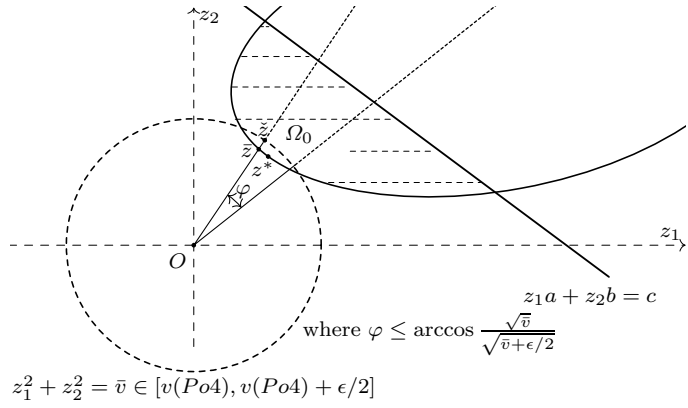

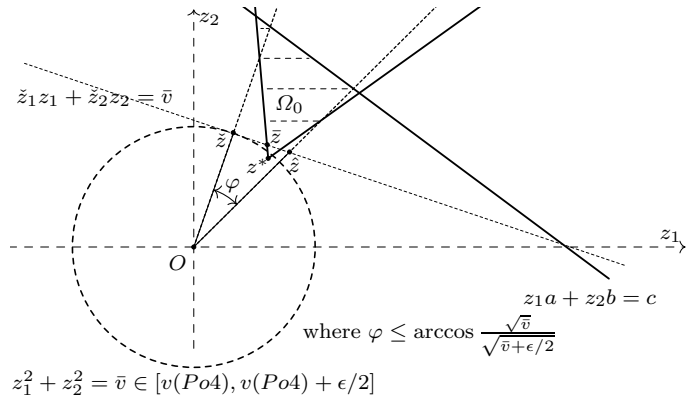
\begin{figure}   %line width=1.pt,dash
	\centering
	\begin{tikzpicture}[line cap=round,line join=round,scale=0.42]\label{fig.1}
	\clip(-5.7458,-4.8) rectangle (15.4,7.548);
	
	\draw [color=black,rotate around={0.:(0.,0.)},line width=0.6pt,dash pattern=on 2.1pt off 2.1pt] (0.,0.) ellipse (3.8cm and 3.8cm);
	
	%	\draw [color=black,rotate around={10.445622:(9.4895238,5.96658)},line width=0.6pt] (9.4895, 5.96658) ellipse (8.4242cm and 4.2379cm);
		
	\draw [->,line width=0.3pt,dash pattern=on 3pt off 3pt] (0,-3.2) -- (0,7.5);
	\draw [->,line width=0.3pt,dash pattern=on 3pt off 3pt] (-5.8,0) -- (15.4,0);
	\draw [color=black, line width=0.7pt] (13,-1) -- (-5,12.5);
		
	\draw [line width=0.2pt,dash pattern=on 1pt off 1pt] (2.1, 6.9)-- (4.98-2.6, 6.9);
	\draw [line width=0.2pt,dash pattern=on 2pt off 2pt] (2.2, 5.957)-- (6.4079-2.8, 5.957);
	\draw [line width=0.2pt,dash pattern=on 2pt off 2pt] (2.3, 4.995)-- (7.673-2.8, 4.995);
	\draw [line width=0.2pt,dash pattern=on 2pt off 2pt] (2.4, 3.9766)-- (4, 3.9766);
	
	%2.33, 2.7999
	
	%C%
	\draw [line width=0.7pt,color=black] (2.33, 2.7999)-- (0, 30);
	\draw [line width=0.7pt,color=black] (2.33, 2.7999)-- (40, 30);

	\draw [line width=0.5pt,color=black,dash pattern=on 1pt off 1pt] (0, 0)-- (10, 10);
	
	\draw [line width=0.5pt,color=black,dash pattern=on 1pt off 1pt] (0,0)-- (10, 29);
	
	\draw [line width=0.2pt,color=black,dash pattern=on 1pt off 1pt,domain=-5.3:13.5] plot(\x,{(14.44-1.225*\x)/3.61});
	
	\begin{scriptsize}
	\draw [fill=black] (2.33, 2.7999) circle (1.8pt); %z^*
	\draw [fill=black] (0, 0) circle (1.8pt);
	
	\draw [fill=black] (3.0, 3.0) circle (1.8pt); %hat z
	
	\draw [fill=black] (2.3, 3.23) circle (1.8pt); %\bar z
	
	\draw [fill=black] (1.235, 3.61) circle (1.8pt); %check z
	\end{scriptsize}
	
	\draw (14.93,0.34) node {${z_1}$};
	\draw (0.5,7.21) node {${z_2}$};
	
	\draw (1.9*1.051,2.45*1.051) node {${z^*}$};
	\draw (0.81*1.051,3.2*1.051) node {$\check{z}$};
	\draw (3.15,2.6) node {$\hat{z}$};	
	\draw (2.56, 3.49) node {$\bar{z}$};
	
	\draw (3,4.5) node {$\Omega_0$};
	\draw (-0.5,-0.5) node {$O$};	
	
	\draw (-3.1,4.7) node[color=black] {$\check{z}_1z_1+\check{z}_2z_2={\bar{v}}$};
	\draw (0,-4.3) node {${z_1^2+ z_2^2=\bar{v}\in[v(Po4), v(Po4)+\epsilon/2]}$};
	\draw (12.4,-1.6) node {${z_1a+z_2b= c}$};
	\draw (7.4,-2.8) node {where $\varphi\leq\arccos\frac{\sqrt{\bar{v}}}{\sqrt{\bar{v}+\epsilon/2}}$};
	
	% danh dau goc
	\draw
	(3,3) coordinate (a) node[right] {}
	-- (0,0) coordinate (b) node[left] {}
	-- (1.25,3.61) coordinate (c) node[above right] {}
	pic["$\varphi$",draw=black,<->,angle eccentricity=1.2,angle radius=0.8cm] {angle=a--b--c};
	%het danh dau goc	
	
	\end{tikzpicture}
	\caption{Graphic representation for the subproblem \eqref{subp4}: $[O,\check{z}]\cap{\Omega_0}=\emptyset$  and $[O,\hat{z}]\cap{\Omega_0}=\emptyset$ .}
\end{figure}

\subsubsection{Analysis for Algorithm \ref{algo2}}

At the first iteration of Algorithm \ref{algo2}, we use the midpoint $\varphi_1=\pi\in[0,2\pi]$ to divide the feasible set $\Omega_0$ of (Po4) into two parts:
\begin{equation*}
\Omega_0\cap \{z\in \mathbb{R}^2\vert~ z_2\ge 0\} \text{ and } \Omega_0\cap \{z\in \mathbb{R}^2\vert~ z_2\le 0\}
\end{equation*}
and then solve
$$v_{1}=\inf\left\{F(z)\left|~z\in\Omega_0,~z_2\ge 0\right.\right\},$$
which is certainly a type of (Po4).

If $v_{1}\le\bar{v},$
we set $\Omega_1=\Omega_0\cap \{z\in\mathbb{R}^2\vert~ z_2\ge 0\}$ and then update
$u_1=\varphi_1=\pi$ with $l_1=l_0$ unchanged. Otherwise,
set $\Omega_1=\Omega_0\cap \{z\in\mathbb{R}^2\vert~ z_2\le 0\}.$ Set $l_1=\varphi_1=\pi$
whereas keeping $u_1=u_0$ as the same.
At the next iteration, $\Omega_1$ is further divided into two parts by either of the midpoints $\varphi_2=\pi/2\in[0, \pi]$ or $\varphi_2=3\pi/2\in[\pi, 2\pi].$
This is equivalent to adding a new cut $z_1\ge 0$ to $\Omega_1$.
We solve
\begin{equation}\label{bisection.2}
\inf\Big\{F(z)|~z\in \Omega_1\cap \{z\in\mathbb{R}^2\vert~ z_1\ge 0\}\Big\}
\end{equation}
to determine whether $(z_1^*,z_2^*)^T$ lies in $\Omega_1\cap \{z\in\mathbb{R}^2\vert z_1\ge 0\}$ or not. If affirmative, we set $\Omega_2=\Omega_1\cap \{z\in\mathbb{R}^2\vert~ z_1\ge 0\}.$ Otherwise,
set $\Omega_2=\Omega_1\cap \{z\in\mathbb{R}^2\vert~ z_1\le 0\}.$
Update the search interval $[l_2,u_2]$ accordingly.
At the third iteration, the set $\Omega_2$ is divided into two parts by the line $(\sin(\varphi_3), -\cos(\varphi_3))z=0$ where the midpoint $\varphi_3$ is either $\pi/4$, or $3\pi/4$, or $5\pi/4$, or $7\pi/4$.

In general, at the $k^{th}$ iteration, Algorithm \ref{algo2} divides $[0, 2\pi]$ into $2^k$ parts. Moreover, {an optimal solution $z^*$ at step 2 of Algorithm \ref{algo1} is} located in $\Omega_k,$ which is a very narrow sector of arc angle $2\pi/2^k$ intersecting with $\Omega_0$. The stopping criterion requires the angle to be smaller than a preset constant, say
$2\pi/2^k\le \arccos\frac{\sqrt{\bar{v}}}{\sqrt{\bar{v}+\epsilon/2}}.$ Therefore the number of steps necessarily to achieve this is
\begin{equation}\label{k-star}
k^*(\bar{v},\epsilon)=\left[\log_2\frac{2\pi}{\arccos\frac{\sqrt{\bar{v}}}
{\sqrt{\bar{v}+\epsilon/2}}}\right]+1,
\end{equation}
where $[t]$ is the integer part of the real number $t$.

When the number of the bisection steps $k$ arrives $k^*(\bar{v},\epsilon),$ we choose the following two points from the boundary of the sector (see Step 2.4):
\begin{eqnarray*}
\check{z}&=&(\sqrt{\bar{v}}\cos u_{k},\sqrt{\bar{v}}\sin u_{k})\in\{z\in \mathbb{R}^2|~ F(z)=z_1^2+z_2^2={\bar v}\},\\
\hat{z}&=&(\sqrt{\bar{v}}\sec{(2\pi/2^k)}\cos u_{k},\sqrt{\bar{v}}\sec{(2\pi/2^k)}\sin u_{k}).
\end{eqnarray*}
Notice that both points $\check z,\hat z$ lie on the line $\{z: \check{z}_1z_1+\check{z}_2z_2={\bar{v}}\},$ which is tangent to the circle $\{z\in \mathbb{R}^2|~ F(z)=z_1^2+z_2^2={\bar v}\}$ at $\check{z}.$ Moreover, for all $\xi\in [\check z,\hat z],$
$$\sqrt{\bar v}=\|\check{z}\|\leq \sqrt{F({\xi})}
\leq \|\hat{z}\| 
= \sqrt{\bar v}\sec\frac{2\pi}{2^k} 
\le \sqrt{\bar v}\cdot\frac{\sqrt{\bar{v}+\epsilon/2}}
{\sqrt{\bar{v}}} 
\le \sqrt{v(\rm{Po4})+\epsilon}.$$
%
%\begin{eqnarray}
%\sqrt{\bar v}=\|\check{z}\|\leq \sqrt{F({\xi})}
%&\leq& \|\hat{z}\| \label{est-0}\\
%&=& \sqrt{\bar v}\sec\frac{2\pi}{2^k} \nonumber\\
%&\le& \sqrt{\bar v}\cdot\frac{\sqrt{\bar{v}+\epsilon/2}}
%{\sqrt{\bar{v}}} \nonumber \\
%&\le& \sqrt{v(\rm{Po4})+\epsilon}. \label{est-1}
%\end{eqnarray}
That is, 
\begin{equation}\label{est-0a}
{F(\xi)}\le v({\rm{Po4}})+\epsilon,~\xi\in [\check z,\hat z].
\end{equation}

Then, depending on the configuration of $\Omega_0,$ the two segments $[O\check{z}]$ and $[O\hat{z}]$ may (Fig. 2) or may not (Fig. 3) intersect $\Omega_0$.
%
%In some extreme cases, for example, when $\Omega_0$ is a very small convex region around $z^*,$ both $[O,\check{z}]\cap\Omega_0$ and $[O,\hat{z}]\cap{\Omega_0}$ can be empty when the algorithm stops. Then, Algorithm \ref{algo2} requires to find a point $\bar{z}$ in $\{z: \check{z}_1z_1+\check{z}_2z_2={\bar{v}}\}\cap \Omega_0$ nearest to $\check{z}$ (see Fig. 2).

- If either $[O,\check{z}]\cap{\Omega_0}\ne\emptyset$ or $[O,\hat{z}]\cap{\Omega_0}\not=\emptyset$, say $[O,\check{z}]\cap{\Omega_0}\ne\emptyset,$ Algorithm \ref{algo2} finds $\bar{z}$ in
$\{z: \check{z}_2z_1-\check{z}_1z_2=0\}\cap \Omega_0$ nearest to $O$. See Fig. 2. Since the line $\{z: \check{z}_2z_1-\check{z}_1z_2=0\}$ passes both $O$ and $\check{z},$ there is
\begin{equation}\label{est-2}
v({\rm{Po4}})=F(z^*)\le F(\bar{z}) \le F({\check z})={\bar v} \le v(\rm{Po4})+\frac{\epsilon}{2}.
\end{equation}
Similarly, if $[O,\hat{z}]\cap{\Omega_0}\not=\emptyset$, Algorithm \ref{algo2} finds $\bar{z}$ such that, by \eqref{est-0a},
\begin{equation}\label{est-3}
v({\rm{Po4}})=F(z^*)\le F(\bar{z}) \le F({\hat z})\le v(\rm{Po4})+{\epsilon}.
\end{equation}
In any case, we have
\begin{equation}\label{est-4}
v({\rm{Po4}})\le F(\bar z)\le v(\rm{Po4})+{\epsilon}.
\end{equation}

- On the other hand, if both $[O,\check{z}]\cap\Omega_0$ and $[O,\hat{z}]\cap{\Omega_0}$ are empty, Algorithm \ref{algo2} requires to find a point $\bar{z}$ in $\{z: \check{z}_1z_1+\check{z}_2z_2={\bar{v}}\}\cap \Omega_0$ nearest to $\check{z}$ (see Fig. 3). The following proposition shows that, when $\Omega_0$ is either unbounded or contains an interior point, such a nearest point $\bar z$ always exists.

\begin{proposition}\label{rezzz22} Let $\check{z}, \hat{z}$ be defined as in Step 2.4 of Algorithm \ref{algo2}. Suppose, at Step 2.5 of Algorithm \ref{algo2}, there happen that 
$[O,\check{z}]\cap{\Omega_0}=\emptyset$ and $[O,\hat{z}]\cap{\Omega_0}=\emptyset.$
\begin{itemize}
\item[-] If $\Omega_0$ is unbounded, then $[\check{z},\hat{z}]\cap{\Omega_0}
\ne\emptyset$.
\item[-] If $\Omega_0$ contains an interior point, then $[\check{z},\hat{z}]\cap{\Omega_0}\ne\emptyset$ for all sufficiently small $\epsilon>0$ where $\epsilon$ is defined as in \eqref{vbar}.	
\end{itemize}
\end{proposition}
\proof
In considering the triangle $\Delta(O\check{z}\hat{z}),$ one has $z^*\in {\Omega_0}
\cap\Delta(O\check{z}\hat{z})$ by Algorithm \ref{algo2}. Since $\Omega_0$ is convex and unbounded and the two edges $[O,\check{z}],~[O,\hat{z}]$ of $\Delta(O\check{z}\hat{z})$ do not intersect with 
$\Omega_0$, the edge $[\check{z}, \hat{z}]$ must intersect with 
$\Omega_0$. Namely, $[\check{z},\hat{z}]\cap\Omega_0\ne\emptyset$.

Now, if  $\Omega_0$ contains an interior point $y'\in\Omega_0$ such that the ball $B(y', r)\subset\Omega_0$ for some $r>0.$ Let $0<\epsilon< r^2$ with $\epsilon$ defined as in \eqref{vbar}. From \eqref{est-0a}, we know $\|\check z\|^2\le v({\rm Po4})+\epsilon$ and 
$\|\hat z\|^2\le v({\rm Po4})+\epsilon$	so that the triangle $\Delta(O\check{z}\hat{z})$ is
contained in the ball of center $O$ with radius $\sqrt{v({\rm Po4})+\epsilon}\le\sqrt{v({\rm Po4})}+\sqrt{\epsilon}<\sqrt{v({\rm Po4})}+r.$ On the other hand, for each $y\in B(y', r)\subseteq \Omega_0$, one has $\|y\|\geq \|z^*\|=\sqrt{v({\rm Po4})},$ which implies that 
$\|y'\|\geq \sqrt{v({\rm Po4})}+r.$ Hence, $y'\in\Omega_0$ lies outside the triangle $\Delta(O\check{z}\hat{z})$ whereas $z^*\in \Delta(O\check{z}\hat{z})$. Due to the convexity of $\Omega_0$, the segment $[z^*, y']\subset \Omega_0$. However, our assumption assumes that $[O,\check{z}]\cap{\Omega_0}=\emptyset$ and $[O,\hat{z}]\cap{\Omega_0}=\emptyset$ so that $[\check{z},\hat{z}]\cap[z^*, y']\not=\emptyset$. The result of this proposition follows immediately. \hfill$\Box$

%Theorem \ref{thm4} below summaries the results in \eqref{k-star}-\eqref{est-4} for %analyzing Algorithm \ref{algo2}, which implements Step 2 in Algorithm \ref{algo1}.

{Hence, with mild assumptions, see Proposition \ref{rezzz22}, $\bar{z}$ always exists and $\bar{z}$ belongs to the segment $[\check{z}, \hat{z}]$, so    by \eqref{est-0a} we have $F(\bar{z})\leq ({\rm{Po4}})+\epsilon$.             }

{In summary, we have the following result}.
\begin{theorem}\label{thm4} Suppose $F(z)=z^2_1+z^2_2$ and that Step 1 in Algorithm \ref{algo1} returns an approximate optimal value $\bar{v}\in [v({\rm{Po4}}), v({\rm{Po4}})+\epsilon/2]$ for a sufficiently small $\epsilon>0.$ Assume that the feasible domain of problem {(Po4)}, $\Omega_0$, is either unbounded or contains an interior point. Then, Algorithm \ref{algo2} terminates after $k^*(\bar{v},\epsilon)$ $($defined as in \eqref{k-star}$)$ number of bisection steps with an approximate solution $\bar z\in\Omega_0$ such that $v({\rm{Po4}})\leq F(\bar{z})\leq v({\rm{Po4}})+\epsilon$.
\end{theorem}

%	$\rm (a)$ The necessary number of steps in Algorithm 2 is \[\left[\log_2\frac{2\pi}{\arccos\frac{\sqrt{\bar{v}}}{\sqrt{\bar{v}+\epsilon/2}}}\right]+1.\]
%	
%$\rm (b)$ Let $\bar{x}$ be a solution of equations $\{f(x)=\bar{z}_1, g(x)=\bar{z}_2\}$, where $\bar{z}_1, \bar{z}_2$ are in Algorithm 2. Then  $v(Po4)\leq F(\bar{z})=F(f(\bar{x}), g(\bar{x}))\leq v(Po4)+\epsilon$.

%
%
%\proof (a) See above arguments.
%
%(b) 	
%- If $\check{z}\not\in B(O, \bar{v})\cap \Omega_0$ (see Fig. 2) then by the fact that $\bar{z}$ belongs to the segment $[\check{z}, \hat{z}]$ ($z^*$ belong to a very narrow sector of arc angle $2\pi/2^k$ intersecting with {\bf C}), we have $\min\{\|\check{z}\|, \|\hat{z}\|\}\leq \sqrt{F(\bar{z})}\leq \max\{\|\check{z}\|, \|\hat{z}\|\}$. On the other hand, since $z_1^*z_1+z_2^*z_2=\sqrt{\bar{v}}$ is perpendicular to $O\check{z}$ so $\sqrt{\bar{v}}=\|\check{z}\| \leq \|\hat{z}\|$ and $\|\hat{z}\|\leq\sqrt{\bar{v}}\cos^{-1}(\arccos\frac{\sqrt{\bar{v}}}{\sqrt{\bar{v}+\epsilon/2}})=\frac{\sqrt{\bar{v}}(\sqrt{\bar{v}+\epsilon/2})}{\sqrt{\bar{v}}}=\sqrt{\bar{v}+\epsilon/2}\leq \sqrt{v(Po4)+\epsilon}$ (see \eqref{vbar}). Hence $v(Po4)\leq F(\bar{z})\leq v(Po4)+\epsilon$.
%
%	- If $\check{z}\in B(O, \bar{v})\cap \Omega_0$ (see Fig. 3) then it is trivial to see that $v(Po4)\leq F(\bar{z})\leq F(\check{z})=\bar{v}\in[v(Po4), v(Po4)+\epsilon/2]$, so $v(Po4)\leq F(\bar{z})\leq v(Po4)+\epsilon$.  \hfill$\Box$

\subsubsection{Discussion on Step 2.5 of Algorithm \ref{algo2}}

To determine whether  $\check{z}=(\sqrt{\bar{v}}\cos u_{k},\sqrt{\bar{v}}\sin u_{k})$ (the same for $\hat{z}$) obtained in Step 2.4 of Algorithm \ref{algo2} satisfies the condition $[O\check{z}]\cap{\Omega_0}\ne\emptyset$, we observe that (see Fig. 2 and Fig. 3)
\begin{eqnarray}
[O\check{z}]\cap{\Omega_0}=\emptyset &\iff& {\inf}\{\check{z}_1z_1+\check{z}_2z_2: \check{z}_2z_1-\check{z}_1z_2=0,~z\in\Omega_0\}>{\bar{v}} \nonumber\\
&\iff& {\bar{v}}<\left\{
\begin{array}{lll}\label{Step2.5-i}
&\inf\limits_{x\in \mathbb{R}^n} & \check{z}_1f(x)+\check{z}_2g(x)  \\
& ~~~~~{\rm s.t.} & f(x) a + g(x) b -c\leq 0; \\
& & \check{z}_2f(x)-\check{z}_1g(x)=0.
\end{array}\right.
\end{eqnarray}
%here we default ${\inf}\{\check{z}_1z_1+\check{z}_2z_2$ on a empty set is unbounded.
Problem \eqref{Step2.5-i} is again a type of (Po4) so that whether or not $\check{z}\in \Omega_0$ can be checked by computing the optimal value $v(\eqref{Step2.5-i})$ of a related SDP as specified in Theorem \ref{thm3asdf}. However, if $\check{z}\not\in \Omega_0,$ finding a point $\bar{z}$ in $\{z: \check{z}_1z_1+\check{z}_2z_2={\bar{v}}\}\cap \Omega_0$ nearest to $\check{z}$ is generally difficult. It amounts to solving a point $\bar x\in\mathbb{R}^n$ such that $\bar{z}_1=f(\bar x),~\bar{z}_2=g(\bar x)$ and $\bar z=(\bar{z}_1,\bar{z}_2)$ is a minimum solution to
\begin{equation}\label{Step2.5-j}
\inf\{\check{z}_2z_1-\check{z}_1z_2: \check{z}_1z_1+\check{z}_2z_2={\bar{v}},~z\in\Omega_0\}.
\end{equation}

Nevertheless, there is one special case of \eqref{Step2.5-j} which we know how to circumvent the difficulty.
When $a=b=c=0,$ \eqref{Step2.5-j} becomes
\begin{equation}\label{Step2.5-k}
\inf\{\check{z}_2f(x)-\check{z}_1g(x): \check{z}_1f(x)+\check{z}_2g(x)={\bar{v}},~x\in\mathbb{R}^n\}.
\end{equation}
It is a (QP1EQC). See Introduction. By applying the S-lemma with equality \cite{XWS} and the standard rank-one decomposition procedure \cite{sz}, we can obtain an optimal solution $\bar x$ to \eqref{Step2.5-k} such that $\bar z=(f(\bar x),g(\bar x))\in\Omega_0={\bf C}$ is nearest to $\check{z}$. According to Theorem \ref{thm4}, there indeed is
$v({\rm{Po4}})\leq F(f(\bar x),g(\bar x))\leq v({\rm{Po4}})+\epsilon.$

Notice that, in this special case $a=b=c=0,$ Step 3 in Algorithm \ref{algo1} for solving a simultaneous quadratic system is not necessary. In the following section, we shall see that the special case suffices to solve the
quadratic surfaces intersection problem (QSIC). We thus summarize the discussion into the following theorem.
%
%
%
%\noindent{\bf - Case 2:} The constraint $z_1a+z_2b\leq c$ is not available ((QSIC) is a special case).
%
%Because of the fact that $\bar{z}$ is a boundary point of ${\bf C}$, our algorithm is implementable. %(see Algorithm 2, Theorem \ref{thm4} and \ref{thm5n})
%% at $\bar{z}$, namely $\bar{z}_1a+\bar{z}_2b<c$, where $\bar{z}$ is defined in Algorithm \ref{algo2}.
%Indeed, an approximation of the optimal solution of (Po4) can be computed via  Algorithm \ref{algo2} (see Theorem \ref{thm4}) and the following result (this result ensures that equations $\left\{\begin{array}{ll}
%f(x)=\bar{z}_1\\
%g(x)=\bar{z}_2
%\end{array}\right.$ can be computed  by solving an (SDP)).

\begin{theorem}\label{thm5n} Let $a=b=c=0;$ $\bar v$ be an approximate value of $v({\rm Po4})$ as specified in \eqref{vbar} and $\check{z},~\hat{z}$ be obtained in Step 2.4 of Algorithm \ref{algo2}.\\
\noindent$\bullet$ If $[O\check{z}]\cap{\Omega_0}=\emptyset$ and $[O\hat{z}]\cap{\Omega_0}=\emptyset$, we solve $\bar{x}\in{\rm argmin}\{\check{z}_2f(x)-\check{z}_1g(x): \check{z}_1f(x)+\check{z}_2g(x)={\bar{v}}\};$\\
\noindent$\bullet$ If $[O\check{z}]\cap{\Omega_0}\ne\emptyset$, we solve $\bar{x}\in{\rm argmin}\{\check{z}_1f(x)+\check{z}_2g(x): \check{z}_2f(x)-\check{z}_1g(x)=0\};$ \\
\noindent$\bullet$ If $[O\hat{z}]\cap{\Omega_0}\ne\emptyset$, we solve $\bar{x}\in{\rm argmin}\{\hat{z}_1f(x)+\hat{z}_2g(x): \hat{z}_2f(x)-\hat{z}_1g(x)=0\};$\\
where $[O{z}]\cap{\Omega_0}=\emptyset \iff {\bar{v}}<{\inf}\{{z}_1f(x)+{z}_2g(x): {z}_2f(x)-{z}_1g(x)=0,~x\in\mathbb{R}^n\}.$ Then, $\bar x$ is an approximate optimal solution to {\rm(Po4)} with $F(z)=z^2_1+z^2_2$ such that
$v({\rm{Po4}})\leq F(f(\bar x),g(\bar x))\leq v({\rm{Po4}})+\epsilon.$  
\end{theorem}

\section{Applications}

\subsection{Problem (QSIC)} %\textcolor{red}{(QSIC)}???

The problem (QSIC) defined in \eqref{ap:0100} can be recast in the following form:
\begin{eqnarray*}%\label{ap:0100}
		{\rm (QSIC)}\hspace*{0.3cm}
		\inf\limits_{x\in \mathbb{R}^{n}} (x^TPx+p^Tx+p_0)^2+(x^TQx+q^Tx+q_0)^2.
\end{eqnarray*}
It is to determine whether or not the two hypersurfaces, $f(x)=x^TPx+p^Tx+p_0=0$ and
$g(x)=x^TQx+q^Tx+q_0=0$, intersect. If $v{\rm (QSIC)}=0,$ or $v{\rm (QSIC)}<\rho$ for some user-defined small $\rho>0$ to accommodate possible noises arising from the problem data, the two hypersurfaces are thought to intersect with each other.
Otherwise, $v{\rm (QSIC)}$ gives some sort of measurement as to how far the two hypersurfaces deviate from each other. For this purpose, the optimal value $v{\rm (QSIC)}$ suffices to decide the intersection problem. However, problem (QSIC) belongs to the special case that $F(z)=z_1^2+z_2^2$ and $a=b=c=0.$ By Theorem \ref{thm5n}, we are able to find an approximate solution $\bar x$ which is either a common root of the two quadratic equations $f(x)=g(x)=0;$ or can be viewed as a point whose total distance to the two hypersurfaces $f(x)=0$ and $g(x)=0$ is the least.

%In this application, only the optimal value $v{\rm (QSIC)}$ is concerned. Computing an optimal solution to decide the intersection problem is not an option.

We solve ${\rm (QSIC)}$ by two cases.

\noindent{\bf Case 1}: If $\{P, Q\}$ are linearly independent, by Theorem \ref{th001} and by taking
$F(z)=z_1^2+z_2^2$, $a=b=c=0,$
we see that (QSIC) is a (Po4) satisfying conditions \eqref{cond}-\eqref{cond0}.    In this case, the optimal value $v(\rm{QSIC})$
can be computed by the following (SDP):
\begin{align*}\gamma^*&=\underset{\gamma, \,\alpha, \,\beta \in \mathbb{R}}{\sup}  \left\lbrace\begin{array}{cc} \gamma&|~M\succeq 0 \end{array}\right\rbrace,
\end{align*} where $M$ is $$\begin{pmatrix}\begin{array}{ccc}
1 & \quad \quad \quad& 0 \\ 0 &\qquad & 1
\end{array} & [0] & \begin{array}{cc}
\frac{-\alpha}{2}  \\ \frac{-\beta}{2}
\end{array}\\  [0]^T & \alpha P+\beta Q & \alpha p+\beta q\\ \begin{array}{cc}
\frac{-\alpha}{2}  &\quad \quad \frac{-\beta}{2}
\end{array} & \quad \alpha p^T+\beta q^T &\quad \alpha p_0+\beta q_0-\gamma\end{pmatrix}.$$
The optimal solution of (QSIC), if attainable, can be computed approximately via Algorithm \ref{algo1},  Algorithm \ref{algo2} and Theorem \ref{thm5n}.

\vskip 0.2cm

\noindent{\bf Case 2}: If $\{P, Q\}$ are linearly dependent, say $Q=t^*P.$ If $P=Q=0,$ then (QSIC) is an unconstrained convex quadratic optimization problem, which can be solved directly. Hence, let us assume that $P\not=0.$ Then, we multiply the first equation in \eqref{nls:001} by $t^*$ and subtract it from the second equation to obtain
\begin{eqnarray}
\begin{cases}\begin{array}{lll}
x^TPx+p^Tx+p_0-z_1=0,\\
(q^T-t^*p^T)x+(q_0-t^*p_0)+t^*z_1-z_2=0.\label{27uhhu}
\end{array}\end{cases}
\end{eqnarray}

Define $y^T=[x^T,~{z^T}]$; $h^T=[q^T-t^*p^T, t^*, -1]$; $h_0=q_0-t^*p_0$;
and
$$\bar{A}=\begin{pmatrix}\begin{array}{ll}[0]_{n\times n} & [0]_{n\times 2} \\
[0]_{2\times n} &I_{2\times 2}\end{array}\end{pmatrix},~ \bar{P}=\begin{pmatrix}\begin{array}{ll}P & [0]_{n\times 2} \\
[0]_{2\times n} &[0]_{2\times 2}\end{array}\end{pmatrix}\not=0.$$
Then, ${\rm(QSIC)}$ becomes
\begin{eqnarray}
&\inf_{y\in \mathbb{R}^{n+2}}& y^T\bar{A}y  \nonumber\\
&{\rm s.t.}&y^T\bar{P}y+\bar{p}^Ty+p_0=0,~h^Ty+h_0=0,\nonumber
\end{eqnarray}
where $\bar{p}^T=[p^T, -1, 0].$
By the null space representation, the solution set for the hyperplane $h^Ty+h_0=0$
can be written as $y=y_0+Vz$ for some $z\in \mathbb{R}^{n+1},$ where
$y_{0} = -\frac{h_0}{h^{T}h}h$ and $V\in\mathbb{R }^{(n+2) \times (n +1)}$ is the matrix basis of  $\mathcal {N}(h)$.
It implies that ${\rm(QSIC)}$ is reduced to the following (QP1EQC) problem:
\begin{eqnarray}
&\inf\limits_{z\in \mathbb{R}^{n+1}}&~ (y_0+Vz)^T\bar{A}(y_0+Vz)\hskip 1.5cm \label{dfgh}\\
&{\rm s.t.}&~(y_0+Vz)^T\bar{P}(y_0+Vz)+\bar{p}^T(y_0+Vz)+p_0=0,\nonumber
\end{eqnarray}
which can be solved by applying the S-lemma with equality. See \cite{XWS}.

\subsection{Problem (AQP)}
Notice that (AQP) can be equivalently formulated as follows:
%\begin{eqnarray}
%&\inf\limits_{x\in \mathbb{R}^n}& (x^TPx+p^Tx+p_0)^2.  \label{ap:100}\\
%&{\rm s.t.}& x^TQx+q^Tx+q_0\le 0\nonumber
%\end{eqnarray}
%It can be further written as
\begin{eqnarray}\label{ap:100}
&\inf\limits_{(x,z_1)\in \mathbb{R}^n\times\mathbb{R}}& z_1^2 \hskip2.5cm\nonumber\\
&{\rm s.t.}&\hskip-0.5cm\begin{cases}\begin{array}{ll}
x^TPx+p^Tx+p_0-z_1=0, \\
x^TQx+q^Tx+q_0\leq 0. %\label{ap:1000x}
\end{array}
\end{cases}
\end{eqnarray}
We again split the discussion into two cases: Case 1: $\{P, Q\}$ are linearly independent; Case 2: $\{P, Q\}$ are linearly dependent.

\noindent $\bullet$ {\bf Case 1:} If $\{P, Q\}$ are linearly independent, by writing \eqref{ap:100} as
\begin{eqnarray}\label{ap:1000x}
&\inf\limits_{(x,z_1)\in \mathbb{R}^n\times\mathbb{R}}& z_1^2 \hskip2.5cm\nonumber\\
&{\rm s.t.}&\hskip-0.5cm\begin{cases}\begin{array}{ll}
z_2\leq 0,  \\
f(x)=x^TPx+p^Tx+p_0=z_1, \\
g(x)=x^TQx+q^Tx+q_0=z_2,
\end{array}
\end{cases}
\end{eqnarray}
with $F(z)=z_1^2, a=0, b=1, c=0,$  we see that (AQP) is a (Po4) satisfying conditions \eqref{cond}-\eqref{cond0}. Then, the optimal value $v(\eqref{ap:1000x})$ can be computed by:
\begin{align*}\gamma^*=v{\rm (AQP)}&=\underset{\gamma, \,\alpha, \,\beta \in \mathbb{R}}{\sup}  \left\lbrace\begin{array}{cc} \gamma&|~M\succeq 0 \end{array}\right\rbrace,
\end{align*} where $M$ is $$\begin{pmatrix}\begin{array}{ccc}
1 & \quad \quad \quad& 0 \\ 0 &\qquad & 0
\end{array} & [0] & \begin{array}{cc}
\frac{-\alpha}{2}  \\ \frac{\mu-\beta}{2}
\end{array}\\  [0]^T & \alpha P+\beta Q & \alpha p+\beta q\\ \begin{array}{cc}
\frac{-\alpha}{2}  & \quad\quad \frac{\mu-\beta}{2}
\end{array} &\quad \alpha p^T+\beta q^T & \quad\alpha p_0+\beta q_0-\gamma\end{pmatrix}.$$

Although the optimal solution can be computed approximately via Algorithm \ref{algo1} and Algorithm \ref{algo2}, yet due to $F(z)=z_1^2$ and $z_2\leq 0$ being simple functions, we can solve it in a better way without resorting to the bisection method.
\begin{figure}   %line width=1.pt,dash
	\centering
	\begin{tikzpicture}[line cap=round,line join=round,scale=0.38]
	\clip(-5.7458, 0.0) rectangle (19, 8);
	
	\draw [line width=1pt,dash pattern=on 3pt off 3pt] (1.1 , 1.5) -- (1.1 ,7.5);
	\draw [line width=1pt,dash pattern=on 3pt off 3pt] (-1.1 , 1.5) -- (-1.1 ,7.5);
	\draw [color=black,rotate around={10.445622:(9.4895238,5.96658)},line width=1.pt] (9.4895, 5.96658) ellipse (8.4242cm and 4.2379cm);
	\draw [color=black, line width=1pt] (-5, 6)-- (17.8, 6);
	
	\draw [line width=0.2pt,dash pattern=on 2pt off 2pt] (4.1, 5)-- (17, 5);
	\draw [line width=0.2pt,dash pattern=on 2pt off 2pt] (1.3, 3.9766)-- (15.7, 3.9766);
	\draw [line width=0.2pt,dash pattern=on 2pt off 2pt] (2.15, 2.9669)-- (13.8, 2.9669);
	\draw [line width=0.2pt,dash pattern=on 2pt off 2pt] (4.22, 1.98)-- (11, 1.98);
	
	\draw [->,line width=0.3pt,dash pattern=on 3pt off 3pt] (0, 1.5) -- (0,7.5);
	\draw [->,line width=0.3pt,dash pattern=on 3pt off 3pt] (-4.5,6) -- (18.3, 6);
	
	\draw (18.2, 5.5) node {${z_1}$};
	\draw (0,7.8) node {${z_2}$};
	\draw (-0.5,6.5) node {$O$};	
	\draw (0, 1.0) node {${z_1^2=v(\rm{Po4})}$};
	\draw (8, 4.5) node {$\bf{C}$};
	\draw (2.7, 5) node {${(z_1^*, z_2^*)}$};
	\draw (8, 6.4) node {${z_2= 0}$};
	\draw (-4.3, 0.7) node {\bf{Fig}. 4};
	
	\begin{scriptsize}
	\draw [color=black] (4,0) circle (0.5pt);
	\draw [fill=black] (1.13, 5) circle (4.1pt);
	\draw [fill=black] (0, 6) circle (3.1pt);
	\end{scriptsize}	
	\end{tikzpicture}

	\begin{tikzpicture}[line cap=round,line join=round,scale=0.38]
	\clip(-5.7458, 0.0) rectangle (19, 8);
	
	\draw [line width=1pt,dash pattern=on 3pt off 3pt] (1.1 , 1.5) -- (1.1 ,7.5);
	\draw [line width=1pt,dash pattern=on 3pt off 3pt] (-1.1 , 1.5) -- (-1.1 ,7.5);	
	\draw [line width=1pt] (8,8) circle (7.2cm);	
	\draw [color=black, line width=1pt] (-5, 6)-- (17.8, 6);	
	
	\draw [line width=0.2pt,dash pattern=on 2pt off 2pt] (1.55, 5)-- (14.6, 5);
	\draw [line width=0.2pt,dash pattern=on 2pt off 2pt] (1.99, 3.9766)-- (14, 3.9766);
	\draw [line width=0.2pt,dash pattern=on 2pt off 2pt] (3, 2.9669)-- (13.1, 2.9669);
	\draw [line width=0.2pt,dash pattern=on 2pt off 2pt] (4.22, 1.98)-- (12, 1.98);
	
	\draw [->,line width=0.3pt,dash pattern=on 3pt off 3pt] (0, 1.5) -- (0,7.5);
	\draw [->,line width=0.3pt,dash pattern=on 3pt off 3pt] (-4.5,6) -- (18.3, 6);
	
	\draw (18.2, 5.5) node {${z_1}$};
	\draw (0,7.8) node {${z_2}$};
	\draw (-0.5,6.5) node {$O$};
	\draw (8, 6.4) node {${z_2= 0}$};
	\draw (0, 1.0) node {${z_1^2=v(\rm{Po4})}$};
	\draw (8, 4.5) node {$\bf{C}$};
	\draw (2.7, 6.6) node {${(z_1^*, z_2^*)}$};	
	\draw (-4.3, 0.7) node {\bf{Fig}. 5};
	
	\begin{scriptsize}
	\draw [color=black] (4,0) circle (0.5pt);
	\draw [fill=black] (1.13, 6) circle (4.1pt);
	\draw [fill=black] (0, 6) circle (3.1pt);
	\end{scriptsize}
	\end{tikzpicture}
\end{figure}

- First, if $\gamma^*=v{\rm (AQP)}=0,$ then any optimal solution to the following
(QP1EQC) is an optimal solution to \eqref{ap:1000x}.
\begin{eqnarray*}%\label{po4}
	{\rm (QP1EQC)}\hspace*{0.3cm}
	\begin{array}{lll}
		&\inf\limits_{x\in \mathbb{R}^n} &g(x)\\
		&~{\rm s.t.}&f(x)= 0.
	\end{array}
\end{eqnarray*}
The solution of (QP1EQC) can be obtained by solving an SDP and a standard matrix rank one decomposition. See \cite{XWS,sz}.

- Secondly, if $\gamma^*=v{\rm (AQP)}>0$, then one has $\{z\in\mathbb{R}^2|~ z_1^2< \gamma^*\}\ne \emptyset$ and
\begin{eqnarray*}
%	 \{z\in\mathbb{R}^2: z_1^2< \gamma^*\}\ne \emptyset \text{ and }\nonumber\\
	 \{z\in\mathbb{R}^2: -\sqrt{\gamma^*}<z_1< \sqrt{\gamma^*}\}\cap \big({\bf C}\cap \{z\in\mathbb{R}^2|~ z_2\leq 0\}\big)=\emptyset. %\label{xcvb}
\end{eqnarray*}

Since ${\bf C}\cap \{z\in\mathbb{R}^2: z_2\leq 0\}$ is convex, there must be either
$$\big({\bf C}\cap \{z\in\mathbb{R}^2|~ z_2\leq 0\}\big) \subset \{z\in\mathbb{R}^2|~ z_1\geq \sqrt{\gamma^*}\};$$
or
$$\big({\bf C}\cap \{z\in\mathbb{R}^2|~ z_2\leq 0\}\big) \subset \{z\in\mathbb{R}^2|~ z_1\leq -\sqrt{\gamma^*}\}.$$
In the former case, (AQP) is solved by $\inf\{f(x): g(x)\le0\}$ (see Fig. 4 and 4), while the latter case can be solved by $\inf\{-f(x): g(x)\le0\}.$ Both problems are (QP1QC).

\noindent $\bullet$ {\bf Case 2:} If $\{P, Q\}$ are linearly dependent such that $Q=t^*P$, we assume that there exists a minimizer where KTCQ is  satisfied. Now, we can use a linear combination to eliminate the matrix $Q$ in \eqref{ap:1000x} and get
\begin{eqnarray}
&\inf\limits_{(x,z_1)\in \mathbb{R}^n\times\mathbb{R}}& z_1^2 \hskip 5.5cm \nonumber \\ %\label{ap:10001}\\
&{\rm s.t. } &  \begin{cases}\begin{array}{ll}
x^TPx+p^Tx+p_0-z_1=0,\\
(q^T-t^*p^T)x+(q_0-t^*p_0)+t^*z_1 \leq 0. \label{ap:10000}\\
\end{array}\end{cases}
\end{eqnarray}
Then, we solve an optimal solution, say $x^*$, to (AQP) (if attainable) by fully analyzing the KKT system of \eqref{ap:10000}. To this end, let us assume some suitable constraint qualification hold at $x^*.$

Let $\lambda_1\in\mathbb{R},~\lambda_2\ge0$ be the Lagrange multipliers associated with the two constraints in \eqref{ap:10000}, respectively. The KKT system can be written down as follows:
\begin{eqnarray}\label{35poiu}\begin{cases}
\begin{array}{l}
-2z_1=-\lambda_1+\lambda_2t^*;\\
0=\lambda_1(2Px+p)+\lambda_2(q^T-t^*p^T);\\
x^TPx+p^Tx+p_0-z_1=0;\\
(q^T-t^*p^T)x+(q_0-t^*p_0)+t^*z_1 \leq 0;\\
\lambda_2\big((q^T-t^*p^T)x+(q_0-t^*p_0)+t^*z_1\big)=0;\\
\lambda_2\geq 0.
\end{array}
\end{cases}\end{eqnarray}

- Case 1: $\lambda_2> 0$. By the complementary slackness, $(q^T-t^*p^T)x+(q_0-t^*p_0)+t^*z_1=0.$ That is, the second inequality constraint in \eqref{ap:10000} is active. Therefore, to check whether there are KKT candidates $(z_1^*,x^*)$ associated with $\lambda_2> 0$ satisfying \eqref{35poiu}, we can instead find optimal solutions $(z_1^*,x^*)$ to the following problem %\eqref{ap:10005}
\begin{eqnarray}\label{ap:10005}
\begin{array}{ll}
\inf\limits_{(x,z_1)\in \mathbb{R}^n\times\mathbb{R}}& z_1^2  \\
  ~~~~~~~{\rm s.t.} & x^TPx+p^Tx+p_0-z_1=0; \\
 & (q^T-t^*p^T)x+(q_0-t^*p_0)+t^*z_1= 0,
\end{array}
\end{eqnarray}
which can be reduced to an (QP1EQC) problem as in \eqref{dfgh}.

%
%\begin{eqnarray}\eqref{ap:100051}
%\begin{array}{lll}
%&\inf\limits_{(z_1, x^T)^T\in \mathbb{R}^{n+1}} & z_1^2  \\
% & ~~~~~{\rm s.t.} & x^TPx+p^Tx+p_0-z_1=0; \\
%& & (q^T-t^*p^T)x+(q_0-t^*p_0)+t^*z_1= 0,
%\end{array}
%\end{eqnarray}

- Case 2: $\lambda_2 = 0$, then (\ref{35poiu}) becomes
\begin{eqnarray}\label{36hgfd}\begin{cases}
\begin{array}{l}
-2z_1=-\lambda_1;\\
0=\lambda_1(2Px+p);\\
x^TPx+p^Tx+p_0-z_1=0;\\
(q^T-t^*p^T)x+(q_0-t^*p_0)+t^*z_1 \leq 0.
\end{array}
\end{cases}\end{eqnarray}

We now further divide this case to two sub-cases:
$\lambda_1 = 0$ and $\lambda_1 \ne 0$.

With $\lambda_1 =\lambda_2 = 0,$ we have $z_1=0$ and (\ref{36hgfd}) becomes
\begin{eqnarray}\label{37hgfd}\begin{cases}
\begin{array}{l}
x^TPx+p^Tx+p_0=z_1=0;\\
-(q^T-t^*p^T)x+(q_0-t^*p_0)\leq 0.
\end{array}
\end{cases}\end{eqnarray}
Since the objective function is $z_1^2$, any solution of (\ref{37hgfd}) is also a solution of \eqref{ap:10000}. It can be decided by solving the following (QP1EQC):
%\begin{eqnarray}
%&\inf\limits_{x\in \mathbb{R}^{n}}&~ -(q^T-t^*p^T)x+(q_0-t^*p_0)\hskip 1.5cm \label{dfgh.01}\\
%&{\rm s.t.}&~x^TPx+p^Tx+p_0=0.\nonumber
%\end{eqnarray}
\begin{eqnarray}\label{dfgh.0}
\begin{array}{lll}
&\inf\limits_{x\in \mathbb{R}^{n}} & -(q^T-t^*p^T)x+(q_0-t^*p_0)  \\
 & ~~{\rm s.t.} & x^TPx+p^Tx+p_0=0.
\end{array}
\end{eqnarray}

If the optimal value $v(\eqref{dfgh.0})>0,$ then
(\ref{37hgfd}) does not have a solution. In this case, $\lambda_1 \ne 0, \lambda_2 = 0$ such that \eqref{36hgfd} becomes
\begin{eqnarray}\label{38hgfd}\begin{cases}
\begin{array}{l}
2z_1=\lambda_1\ne 0;\\
0=\lambda_1(2Px+p);\\
x^TPx+p^Tx+p_0-z_1=0;\\
(q^T-t^*p^T)x+(q_0-t^*p_0)+t^*z_1 \leq 0.
\end{array}
\end{cases}\end{eqnarray}
Then, $2Px+p=0$ and
$$z_1=x^TPx+p^Tx+p_0=\frac{1}{2}p^Tx+p_0.$$
Hence, \eqref{38hgfd} is reduced to the following linear system in variables $(x,z_1,\lambda_1)\in\mathbb{R}^{n+2}$:
\begin{eqnarray}\label{39hgfd}\begin{cases}
\begin{array}{l}
2z_1=\lambda_1\ne 0;\\
2Px+p=0;\\
\frac{1}{2}p^Tx+p_0-z_1=0;\\
(q^T-t^*p^T)x+(q_0-t^*p_0)+t^*z_1 \leq 0.
\end{array}
\end{cases}\end{eqnarray}

In summary, when $\{P, Q\}$ are linearly dependent, we can solve (i) $\lambda_2>0$ and \eqref{ap:10005}; (ii) $\lambda_1=\lambda_2=0$ and \eqref{dfgh.0} with $v(\eqref{dfgh.0})\le0$; and (iii) $\lambda_1\not=0,~\lambda_2=0$ and the linear system \eqref{39hgfd}. The optimal solution of (AQP) will be the valid ones among (i), (ii) and (iii) with the smallest objective value $z_1^2$ of (AQP).

\begin{example}
\begin{eqnarray}\label{khkh000000110}
\begin{array}{ll}
\inf\limits_{x\in \mathbb{R}^{2}} & \big|-x_1^2+x_2^2+x_1\big|  \\
  ~~{\rm s.t.} & -x_1^2+x_2^2+1\leq 0.
\end{array}
\end{eqnarray}
	In this example, $\{P,Q\}$ are linearly dependent with $t^*=1.$ Using the linear combination with $t^*=1$ to eliminate $Q$, we can write
\begin{eqnarray}\label{qwert3}
\begin{array}{ll}
\inf\limits_{ (x,z_1)\in \mathbb{R}^2\times\mathbb{R}}& z_1^2  \\
  ~~~~~~~{\rm s.t.} & -x_1^2+x_2^2+x_1-z_1=0, \\
 & -x_1+z_1+1\leq 0.
\end{array}
\end{eqnarray}
Since we are going to utilize the KKT system of \eqref{qwert3}, we first check the constraint qualification. Notice that, $$\nabla(-x_1^2+x_2^2+x_1-z_1)=(-2x_1+1, 2x_2, -1)^T,~~\nabla(-x_1+z_1+1)=(-1, 0, 1)^T.$$ The two gradients are linearly dependent if and only if $x_1=x_2=0$. When $x_1=x_2=0$, \eqref{qwert3} has an empty feasible set $\{z_1=0,z_1\le-1\}=\emptyset.$ In other words, every optimal solution to problem (\ref{khkh000000110}) must be a regular point of the KKT system.

	Then, check all the following three cases:
	\begin{itemize}
		\item[(i)] $\lambda_2> 0.$ Then, the inequality
		$-x_1+z_1+1\leq 0$ is active so that $x_1=z_1+1.$ We first solve the (QP1EQC):
\begin{eqnarray}\label{238uyyu}
\begin{array}{ll}
\inf\limits_{(z_1,x_2)\in \mathbb{R}\times \mathbb{R}}& z_1^2  \\
  ~~~~~{\rm s.t.} & -(z_1+1)^2+x_2^2+1=0
\end{array}
\end{eqnarray}
to obtain the solution $z^*_1=0, x^*_1=1, x^*_2=0$ to \eqref{238uyyu}.
		From the first equation $-2z_1=-\lambda_1+\lambda_2t^*$ in the KKT system \eqref{35poiu}, since $z^*_1=0,~t^*=1,$ we obtain $$\lambda_1=\lambda_2>0.$$
		From the second equation of \eqref{35poiu}, $(x^*_1,x^*_2)^T$ must satisfy
		$$0=(2Px+p)+(q^T-t^*p^T)$$
		which is not true as $(2Px^*+p)+(q^T-t^*p^T)=(-3,0)^T.$ We conclude that there is no KKT point of \eqref{238uyyu} with $\lambda_2>0.$
		\item[(ii)]$\lambda_1=\lambda_2=0.$ In this case, the KKT system \eqref{35poiu} is reduced to (\ref{37hgfd}), which is
		\begin{eqnarray}\label{37hgfd.e}\begin{cases}
		\begin{array}{l}
		-x_1^2+x_2^2-x_1=z_1=0;\\
		-x_1+z_1+1 \leq 0.
		\end{array}
		\end{cases}\end{eqnarray}
		Since \eqref{37hgfd.e} guarantees that $z_1=0,$ any solution to \eqref{37hgfd.e}, if exists, is an optimal solution to (AQP). It is not difficult to see that \eqref{37hgfd.e} has infinitely many solutions among which we can choose, for example, $z^*_1=0, x^*_1=1, x^*_2=0.$
		\item[(iii)]$\lambda_1\not= 0,~\lambda_2=0.$ Dealing with this case might help to find other solutions to \eqref{khkh000000110}. In this case, we need to check the linear system in (\ref{39hgfd}), which is:
		\begin{eqnarray}\label{39hgfd.e}\begin{cases}
		\begin{array}{l}
		2z_1=\lambda_1\ne 0;\\
		(-2x_1, 2x_2)+(1, 0)=(0,0);\\
		z_1=\frac{x_1}{2};\\
		-x_1+z_1+1 \leq 0.
		\end{array}
		\end{cases}\end{eqnarray}
From the second and the third equation of \eqref{39hgfd.e}, we get $x_1=0.5,~x_2=0,~z_1=0.25.$
	They do not satisfy $-x_1+z_1+1 \leq 0$ though. The system \eqref{39hgfd.e}	does not have a solution, either.
	\end{itemize}
\end{example}

%a special type of quadratically constrained quadratic program (QCQP) which minimizes a quadratic function of quadratic subject to several quadratic constraints.
\section{Conclusion and Discussion}

In this paper, we propose a new type of quadratic optimization problems involving a joint numerical range constraint. There are many natural applications arising from such a formulation. Some applications like the double well potential problem (DWP) \cite{FGLSX,XSFX} were already solved independently, while others like (QSIC) and (AQP) can now be resolved through our new approach. Interestingly, due to the composition of ``quadratic with quadratics,'' the objective function $F(f(x), g(x))$ is indeed a polynomial of degree 4. We hope that our method can be later extended to solve optimization problems involving quartic polynomials.

Notice that, in Theorem \ref{th1}, we can replace $\{f, g\}$ by $k$ quadratic functions $\{f_1, \cdots, f_k\}$. The theoretical results in this paper are still valid (the implementation can be more complicate though), provided the following conditions are satisfied:
\begin{eqnarray*}
\label{cond5.1}&&\{(f_1(x), \cdots, f_k(x))| ~x\in\mathbb{R}^n \} \text{~is convex},\\
\label{cond5.2} && \Theta_{k\times k} \succeq 0,
\end{eqnarray*}
where $\Theta_{k\times k}$ defines a quadratic function $f(z)=z^T\Theta_{k\times k}z+\eta^Tz$ of $k$ variables. It suggests that the convexity of the joint numerical range of $\{f_1, \cdots, f_k\}$ is the central feature for many optimization problems and should be studied carefully in the future.

%
%
%\textcolor{red}{In this paper, we minimize a special type of quadratically constrained quadratic program (QCQP).} The problem involves only two quadratic functions $f$ and $g,$ with their images $z_1=f(x),~z_2=g(x)$ forming a convex quadratic objective function on $\mathbb{R}^2.$  \textcolor{red}{ Inparticular, our method to solve (AQP) is better than available methods solving it in literature}.
%
%\textcolor{red}{Interestingly, problem (Po4) is actually a optimization problem of  quartic polynomial subject to a constraint set defined by $m$ quadratic inequalities:
%\begin{eqnarray*}\label{po4}
%	{\rm (Po4)}\hspace*{0.3cm}
%	\begin{array}{lll}
%		v{\rm (Po4)}=&\inf\limits_{x\in \mathbb{R}^n} &F(f(x), g(x)),\\
%		&{\rm s.t.}&f(x)a+g(x)b-c\leq 0
%	\end{array}
%\end{eqnarray*}
%where $a, b, c\in \mathbb{R}^{m}$ and $f(x)=x^TPx+p^Tx+p_0, g(x)=x^TQx+q^Tx+q_0$ are real-valued quadratic functions
%of $n$-variables with $n\times n$ symmetric matrices $P, Q$ whereas $p,q\in\mathbb{R}^n$ and $p_0,q_0\in \mathbb{R},$
%here, the objective function is a special $4th$ degree polynomial (quartic):
%\begin{eqnarray*}
%	&&F(f(x),g(x))\\
%	&=&\theta_1f(x)^2+2\theta_2f(x)g(x)+\theta_3g(x)^2+\eta_1f(x) +\eta_2g(x)
%\end{eqnarray*}
%with $\theta_1, \theta_2, \theta_3, \eta_1, \eta_2 \in \mathbb{R}.$ }

\section*{Acknowledgements}

Huu-Quang, Nguyen's research work
was sponsored partially by Taiwan MOST 107-2811-M-006-535 and
Ruey-Lin Sheu's research work
was sponsored partially by Taiwan MOST 107-2115-M-006-011-MY2.

\noindent Xia's research was supported by National Natural Science Foundation of China
under grants 11822103, 11571029, 11771056, and Beijing Natural Science Foundation Z180005.

\bibliographystyle{amsplain}

\end{document}